\documentclass[12pt]{article}
\usepackage{amsmath,amssymb,amsthm,amscd,latexsym,}
 \font\goth=eusm10
\usepackage{graphicx}

\makeatletter
\renewcommand{\mod}[1]{\allowbreak \if@display \mkern 8mu \else
\mkern 5mu\fi {\operator@font mod}\,\,#1}
\makeatother

\newcommand{\iso}{\cong}
\newcommand{\ga}{\gamma}
\newcommand{\bc}{\mathbb C}
\newcommand{\bn}{\mathbb N}
 \newcommand{\bq}{\mathbb Q}
\newcommand{\br}{\mathbb R}
\newcommand{\bz}{\mathbb Z}
\newcommand{\bff}{\mathbb F}
\newcommand{\bl}{\mathbb L}
\newcommand{\bk}{\mathbb K}
\newcommand{\bo}{\mathbb O}
\newcommand{\bp}{\mathbb P}
\newcommand{\E}{\mathcal E}
\newcommand{\Oc}{\mathcal O}
\newcommand{\wf}{\widetilde{f}}
\newcommand{\wH}{\widetilde{H}}
\newcommand{\wth}{\widetilde{h}}
\newcommand{\wx}{\widetilde{x}}
\newcommand{\wy}{\widetilde{y}}
\newcommand{\aaa}{\mathbb A}
\newcommand{\ddd}{\mathbb D}
\newcommand{\eee}{\mathbb E}

\DeclareMathOperator{\Aut}{{\rm Aut}\,}
\DeclareMathOperator{\Hom}{Hom}
\DeclareMathOperator{\Mod}{Mod}
\DeclareMathOperator{\Pic}{Pic}
\DeclareMathOperator{\rk}{rk}
\DeclareMathOperator{\ch}{ch}
\DeclareMathOperator{\td}{td}
\DeclareMathOperator{\Tyu}{Tyu}
\DeclareMathOperator{\discr}{{\rm discr}\,}
\DeclareMathOperator{\cha}{{\rm char}\,}
\newcommand\F{\mathcal F}
\newcommand\Hh{\mathcal H}
\newcommand\N{\mathcal N}
\newcommand\La{\mathcal L}
\newcommand\Ka{{\mathcal K}}
\newcommand\G{\mathcal G}
\newcommand\M{\mathcal M}

\newcommand\SSS{\mathfrak S}
\newcommand\AAA{\mathfrak A}
\newcommand\DDD{\mathfrak D}

\begin{document}
\title{Degenerations of K\"ahlerian K3 surfaces with finite
symplectic automorphism groups. III\footnote{With support by Russian
Sientific Fund N 14-50-00005.}}
\date{}
\author{Viacheslav V. Nikulin}
\maketitle

\begin{abstract}
Similarly to our papers \cite{Nik9}, \cite{Nik10}, we classify
degenerations of codimension $2$ and higher of K\"ahlerian K3 surfaces with finite
symplectic automorphism groups. In \cite{Nik9}, \cite{Nik10}, it was done
for the codimension one.
\end{abstract}

\section{Introduction}
\label{sec:introduction}
Similarly to our papers \cite{Nik9}, \cite{Nik10}, we classify
degenerations of codimension 2 and higher
of K\"ahlerian K3 surfaces with finite
symplectic automorphism groups. In \cite{Nik9},
\cite{Nik10}, it was done
for the codimension 1. For readers convinience, in  
Section \ref{sec2.tables}, 
we give the Table 1 from \cite{Nik10}
which gives the main classification of such degenerations
of codimension $1$.

In this variant of the paper, we do classification of degenerations of
codimension $\ge 2$ for
finite symplectic automorphism groups of  order $\ge 6$.

For groups of order $> 8$ and the group $Q_8$ of order $8$, 
it is done in Section \ref{sec2.tables}. It has 
Table 2 which 
contains classification of types and
lattices $S$ of such
degenerations. They are main invariants of degenerations.
Table \ref{table3} of Section \ref{sec2.tables} 
gives markings of such degenerations by
Niemeier lattices $N_j$ in notations of \cite{Nik8}, \cite{Nik9}
and \cite{Nik10}. For such groups, the codimension of degenerations is less or
equal to $4$, and for the most of them it is less or equal to $3$.

We consider remaining difficult groups $D_8$ and $(C_2)^3$ of order $8$ in 
Sections \ref{sec3.tables} and \ref{sec4.tables}.  
There are too many cases for these difficult groups, and it is better to 
consider them separately.  

We hope to consider remaining groups $D_6$, $C_4$, $(C_2)^2$, $C_3$, $C_2$ and 
give more results and details in further
variants of the paper and further publications.

%%%%%%%%%%%%%%%%%%%%%%%%%%%%%%%%%%%%%%%%%%%%%%%%%%%%%%%%%
%%%%%%%%%%%%%%%%%%%%%%%%%%%%%%%%%%%%%%%%%%%%%%%%%%%%%%%%%
%%%%%%%%%%%%%%%%%%%%%%%%%%%%%%%%%%%%%%%%%%%%%%%%%%%%%%%%%

\section{Classification of degenerations of\\
K\"ahlerian K3 surfaces with finite symplectic\\
automorphism groups. The main theory.}
\label{sec:theory}

Let $X$ be a K\"ahlerian K3 surface (e. g. see \cite{Sh}, \cite{PS},
\cite{BR}, \cite{Siu}, \cite{Tod}
about such surfaces). That is $X$ is a non-singular compact
complex surface with the trivial canonical class $K_X$, and
its irregularity $q(X)$ is equal to 0.
Then $H^2(X,\bz)$ with the intersection pairing
is an even unimodular lattice $L_{K3}$
of the signature $(3,19)$. The primitive sublattice
$S_X=H^2(X,\bz)\cap H^{1,1}(X)\subset H^2(X,\bz)$
is the {\it Picard lattice} of $X$ generated by first Chern
classes of all line bundles over $X$.

Let $G$ be a finite symplectic automorphism group of $X$.
Here symplectic means
that for any $g\in G$, for a non-zero holomorphic $2$-form $\omega_X\in
H^{2,0}(X)=\Omega^2[X]=\bc\omega_X$, one has $g^\ast(\omega_X)=\omega_X$.
For an $G$-invariant sublattice $M\subset H^2(X,\bz)$, we denote by
$M^G=\{x\in M\ |\ G(x)=x\}$ the {\it invariant sublattice of $M$,} and by
$M_G=(M^G)^\perp_M$ the {\it coinvariant sublattice of $M$.}
By \cite{Nik-1/2}, the coinvariant lattice $S_G=H^2(X,\bz)_G=(S_X)_G$ is
{\it Leech type lattice:} i. e.  it is negative definite,
it has no elements with square $(-2)$, $G$ acts trivially
on the discriminant group $A_{S_G}=(S_G)^\ast/S_G$, and $(S_G)^G=\{0\}$.
For a general pair $(X,G)$, the $S_G=S_X$, and non-general $(X,G)$
can be considered as K\"ahlerian K3 surfaces
with the condition $S_G\subset S_X$
on the Picard lattice (in terminology of \cite{Nik-1/2}).
The dimension of their
moduli is equal to $20-\rk S_G$.

Let $E\subset X$ be a non-singular irreducible rational curve
(that is $E\cong \bp^1$).
It is equivalent to: $\alpha=cl(E)\in S_X$, $\alpha^2=-2$,
$\alpha$ is effective
and $\alpha$ is numerically effective: $\alpha\cdot D\ge 0$
for every irreducible curve
$D$ on $X$ such that $cl(D)\not=\alpha$.

Let us consider the primitive sublattice $S=[S_G, G(\alpha)]_{pr}\subset S_X$
of $S_X$ generated by the coinvariant sublattice $S_G$ and
all classes of the orbit $G(E)$.
We remind that primitive means that $S_X/S$ has not torsion.
Since $S_G$ has no elements with square $(-2)$, it follows
that $\rk S=\rk S_G+1$
and $S=[S_G,\alpha]_{pr}\subset S_X$.

Let us {\it assume that the lattice $S=[S_G,\alpha]_{pr}$
is negative definite.} Then the
elements $G(\alpha)$ define the basis of the root system $\Delta(S)$ of
all elements with square $(-2)$ of $S$. All curves $G(E)$ of $X$
can be contracted to Du Val singularities of types of
connected components of the Dynkin diagram of the basis. The group $G$
will act on the corresponding singular K3 surface
$\overline{X}$ with these Du Val
singularities. For a general such triplet $(X,G,G(E))$,
the Picard lattice $S_X=S$,
and such triplets can be considered as {\it a degeneration
of codimension $1$} of
K\"ahlerian K3 surfaces $(X,G)$ with the finite symplectic
automorphism group $G$.
Really, the dimension of moduli of K\"ahlerian K3 surfaces
with the condition $S\subset S_X$
on the Picard lattice is equal to $20-\rk S=20-\rk S_G-1$.

By Global Torelli Theorem for K3 surfaces \cite{PS},
\cite{BR}, the main invariants
of the degeneration is the {\it type of the abstract group $G$}
which is equivalent
to the isomorphism class of the coinvariant lattice $S_G$,
and the type of the degeneration which is
equivalent to the Dynkin diagram of the basis $G(\alpha)$
or the Dynkin diagram
of the rational curves $G(E)$.

We can consider only the maximal finite symplectic
automorphism group $G$ with the
same coinvariant lattice $S_G$, that is $G=Clos(G)$.
By Global Torelly Theorem for K3 surfaces, this
is equivalent to
$$
G|S_G=\{ g\in O(S_G)\ |\ g\ is\ identity\ on\ A_{S_G}=(S_G)^\ast/S_G \}.
$$
Indeed, $G$ and $Clos(G)$ have the same lattice $S_G$, the same orbits
$G(E)$ and $Clos(G)(E)$, and the same sublattice $S\subset S_X$.

In \cite{Nik9}, all types of $G=Clos(G)$ and types of degenerations
(that is Dynkin
diagrams of the orbits $G(E)$) are described. They are described
in Table 1 below which
is the same as \cite[Table 1]{Nik10}
where ${\bf n}$ gives types of possible $G=Clos(G)$, and we show all possible
types of degenerations at the corresponding rows by their Dynkin diagrams.

In Table 1, the type of $G=Clos(G) $ and the isomorphism class of the lattice $S_G$
is marked by ${\bf n}$. We also give the genus of $S_G$ which
is defined by
the discriminant quadratic form $q_{S_G}$. They
were calculated in papers \cite{Nik0},
\cite{Muk}, \cite{Xiao}, \cite{Hash}.

In \cite{Nik10}, with two exceptions which are marked by $I$ and $II$
for ${\bf n}=10$ (group $D_8$) and the degeneration $2\aaa_1$,
and for ${\bf n}=34$ (group $\SSS_4$) and the degeneration $6\aaa_1$, it was
shown that the lattice $S=[S_G,\alpha]_{pr}$ is unique,
up to isomorphisms, and its
genus (equivalent to $\rk S$
and the discriminant quadratic form $q_S$) is shown in Table 1.
In \cite[Table 2]{Nik10},
we described possible markings of $S$, $G$ and $\alpha$ by
Niemeier lattices which
give exact lattices descriptions of $S$, and the action of
$G$ on $S$ and $G(\alpha)$.

For degenerations of arbitrary codimension $t\ge 1$ which
we consider in this paper,
instead of one orbit $G(E)$
of a non-singular rational curve, we should consider
$t\ge 1$ different orbits
$G(E_1),\dots, G(E_k)$ of non-singular rational curves on $X$, and their
classes $G(\alpha_1), \dots, G(\alpha_k)$ in $S_X$, but we also assume that
the sublattice $S=[S_G,G(\alpha_1),\dots,G(\alpha_k)]_{pr}\subset S_X$
is negative definite. Then the codimension of the degeneration is
equal to $t=\rk S-\rk S_G$, and $S_X=S=[S_G,\alpha_1,\dots, \alpha_t]_{pr}$
in general. We remark that $\rk S=\rk S_G+t\le 19$ since $H^2(X,\bz)$ has
the signature $(3,19)$. Thus, $t\le 19-\rk S_G$.

The type of the degeneration is given by the Dynkin
diagrams and subdiagrams
$$
(Dyn(G(\alpha_1)),\dots, Dyn(G(\alpha_t)))\subset
Dyn(G(\alpha_1)\cup\dots \cup G(\alpha_t))
$$
and their types. In difficult cases, we also
consider the matrix of subdiagrams
which is defined by
$$
(Dyn (G(\alpha_i)),\, Dyn ( G(\alpha_j)))\subset Dyn(G(\alpha_i)\cup G(\alpha_j))
$$
and their types for $1\le i<j\le t$.

In Table 2 which is similar to Table 1 for codimension 1,
we give classification of types of degenerations of arbitrary
codimension $\ge 2$ for ${\bf n}\ge 12$ (then the
codimension is less or equal to $4$).
Equivalently, either $|G|>8$ or $G\cong Q_8$ of order $8$ ($G$ is big enough).
We calculate
genuses of the lattices $S$ by $\rk S$ and $q_S$. By $\ast$,
we mark cases when we prove that
the lattice $S$ is unique up to isomorphisms, for the given type.
In Table 3 which is similar to \cite[Table 2]{Nik10} for the codimension 1,
we give the description of their markings
$S\subset N_j$ by
Niemeier lattices $N_j$, $j=1,\dots, 24$, where we use
results and notations of \cite{Nik7}, \cite{Nik8} and \cite{Nik9}. For 
remaining groups $D_8$ and $(C_2)^3$ of order $8$, similar results are 
given in Sections \ref{sec3.tables} and \ref{sec4.tables}. 

All calculations are similar to our calculations in \cite{Nik9} and
\cite{Nik10} for codimension 1, and they use
similar to \cite{Nik9} and \cite{Nik10}
computer programs.

We hope to give more details and other results in further variants of the paper
and further publications.

\vskip1cm

\section{Types of degenerations for $\bf n\ge 12$}
\label{sec2.tables}

We remind that Niemeier lattices are negative definite
even unimodular lattices
of the rank $24$. There are $24$ such lattices, up to isomorphisms,
$N=N_j$, $j=1,2,\dots,24$, classified by Niemeier. They are characterized
by their root sublattices $N^{(2)}$ generated by all their elements
with square $(-2)$ (called roots). Further, $\Delta(N)$ is the
set of all roots of $N$.
We have the following list of Niemeier lattices $N_j$ where the number $j$
is shown in the bracket:
$$
N^{(2)}=[\Delta(N)]=
$$
$$
(1)\ D_{24},\ (2)\ D_{16}\oplus E_8,\ (3)\ 3E_8,\ (4)\ A_{24},\
(5)\ 2D_{12},\ (6)\ A_{17}\oplus E_7, \ (7)\ D_{10}\oplus 2E_7,
$$
$$
(8)\ A_{15}\oplus D_9,\
(9)\ 3D_8,\
(10)\ 2A_{12},\ (11)\ A_{11}\oplus D_7\oplus E_6,\ (12)\ 4E_6,\
(13)\ 2A_9\oplus D_6,\
$$
$$
(14)\ 4D_6,\
(15)\ 3A_8,\ (16)\ 2A_7\oplus 2D_5,\ (17)\ 4A_6,\ (18)\ 4A_5\oplus D_4,\
(19)\ 6D_4,
$$
$$
(20)\ 6A_4,\
(21)\ 8A_3,\ (22)\ 12A_2,\ (23)\ 24A_1
$$
give $23$ Niemeier lattices $N_j$. The last is Leech lattice (24)
with $N^{(2)}=\{0\}$ which has no roots.

%%%%%%%%%%%%%%%%%%%%%%%%%%%%%%%%%%%%%%%%%%%%%%%%%%%%%%%%%%%%%%%%%%%
%%%%%%%%%%%%%%%%%%%%%%%%%%%%%%%%%%%%%%%%%%%%%%%%%%%%%%%%%%%%%%%%%%%%%
%%%%%%%%%%%%%%%%%%%%%%%%%%%%%%%%%%%%%%%%%%%%%%%%%%%%%%%%%%%%%%%%%%%%%

%%%%%%%%%%%%%%%%%%%%%%%%%%%%%%%%%%%%%%%%%%%%%%%%%
%%%%%%%%%%%%%%%%%%%%%%%%%%%%%%%%%%%%%%%%%%%%%%%%%
%%%%%%%%%%%%%%%%%%%%%%%%%%%%%%%%%%%%%%%%%%%%%%%%%

\begin{table}
\label{table1}
\caption{Types and lattices $S$ of degenerations of
codimension $1$ of K\"ahlerian K3 surfaces
with finite symplectic automorphism groups $G=Clos(G)$.}

%%\index{Table 1}

%%%\addtocontents{toc}{\contentsline {section}{\tocsection {}{T.1}{Table 1}}
%%%{\pageref{table1}}}

\begin{tabular}{|c||c|c|c|c|c|c|c|c|}
\hline
 {\bf n}& $|G|$& $i$&  $G$  & $\rk S_G$ &$q_{S_G}$  &$Deg$& $\rk S$ &$q_S$ \\
\hline
\hline
 $1$  &$2$ & $1$& $C_2$ &   $8$     &$2_{II}^{+8}$ & $\aaa_1$   &  $9$  &
$2_7^{+9}$\\
\hline
      &    &    &    &              &              & $2\aaa_1$  &  $9$  &
$2_{II}^{-6},4_3^{-1}$\\
\hline
\hline
 $2$  & $3$& $1$&$C_3$  &    $12$   & $3^{+6}$       & $\aaa_1$   &  $13$ &
$2_3^{-1},3^{+6}$\\
\hline
      &    &    &       &           &                & $3\aaa_1$  & $13$  &
$2_1^{+1},3^{-5}$\\
\hline
\hline
 $3$  &$4$ & $2$&$C_2^2$&    $12$   & $2_{II}^{-6},4_{II}^{-2}$& $\aaa_1$   & $13$  &
$2_3^{+7},4_{II}^{+2}$\\
\hline
      &    &    &       &           &                            & $2\aaa_1$& $13$  &
$2_{II}^{-4},4_7^{-3}$\\
\hline
      &    &    &       &           &                & $4\aaa_1$  & $13$  &
$2_{II}^{-6},8_3^{-1}$\\
\hline
\end{tabular}
\end{table}

\begin{table}

\begin{tabular}{|c||c|c|c|c|c|c|c|c|}
\hline
 {\bf n}& $|G|$& $i$&  $G$   & $\rk S_G$ &$q_{S_G}$ &$Deg$& $\rk S$ &$q_S$ \\
\hline
\hline
 $17$&$12$ & $3$&$\AAA_4$&  $16$& $2_{II}^{-2},4_{II}^{-2},3^{+2}$ & $\aaa_1$ &$17$  &
$2_7^{-3},4_{II}^{+2},3^{+2}$\\
\hline
     &     &    &        &           &     &$3\aaa_1$ &  $17$  &
$2_1^{-3},4_{II}^{+2},3^{-1}$\\
\hline
     &     &    &        &           &            &$4\aaa_1$ &  $17$  &
$2_{II}^{-2},8_3^{-1},3^{+2}$\\
\hline
     &     &    &        &           &             &$6\aaa_1$ &  $17$  &
$4_1^{-3},3^{+1}$\\
\hline
     &     &    &        &           &             &$12\aaa_1$&  $17$  &
$2_{II}^{-2},8_1^{+1},3^{-1}$\\
\hline
\hline
 $18$&$12$ &$4$ &$D_{12}$&  $16$&$2_{II}^{+4},3^{+4}$ &  $\aaa_1$&  $17$&
$2_7^{+5},3^{+4}$\\
\hline
     &     &    &        &           &   &  $2\aaa_1$    &  $17$  &
$2_{II}^{+2},4_7^{+1},3^{+4}$\\
\hline

\end{tabular}

\end{table}

\begin{table}

\begin{tabular}{|c||c|c|c|c|c|c|c|c|}
\hline
 {\bf n}& $|G|$& $i$&  $G$   & $\rk S_G$ &$q_{S_G}$       &$Deg$& $\rk S$ &$q_S$    \\
\hline
\hline
 $39$&$32$ &$27$&$2^4C_2$&  $17$     &$2_{II}^{+2},4_0^{+2},8_7^{+1}$ & $4\aaa_1$&$18$&
$4_6^{+4}$\\
\hline
     &     &    &        &           &                                & $8\aaa_1$&  $18$&
$2_{II}^{+2},4_7^{+1},8_7^{+1}$\\
\hline
     &     &    &        &           &                                & $16\aaa_1$&$18$&
$2_{II}^{+2},4_6^{+2}$\\
\hline
\hline
$40$ &$32$ &$49$&$Q_8*Q_8$& $17$     & $4_7^{+5}$ & $8\aaa_1$  & $18$ &
$4_6^{+4}$\\
\hline
\hline
$46$ &$36$ &$9$ &$3^2C_4$ & $18$     & $2_6^{-2},3^{+2},9^{-1}$& $6\aaa_1$   &  $19$  &
$4_7^{+1},3^{+1},9^{-1}$               \\
\hline
     &     &    &        &           &                           & $9\aaa_1$ &  $19$  &
$2_5^{-3},3^{+2}$ \\
\hline
     &     &    &        &           &                           & $9\aaa_2$ &  $19$  &
$2_5^{-1},3^{+2}$\\
\hline

\end{tabular}

\end{table}

\begin{table}
\label{table2}
\caption{Types and lattices $S$ of degenerations of
codimension $\ge 2$ of K\"ahlerian K3 surfaces
with finite symplectic automorphism groups
$G=Clos(G)$ for ${\bf n}\ge 12$.}

\begin{tabular}{|c||c|c|c|c|c|c|c|c|}
\hline
 {\bf n}& $|G|$& $i$&  $G$   & $\rk S_G$ &$Deg$& $\rk S$ &$q_S$            \\
\hline
\hline
  $12$&$8$& $4$ & $Q_8$ &  $17$     &
  $(8\aaa_1,8\aaa_1)\subset 16\aaa_1$ & $19$  &
$2_{II}^{-2},8_5^{-1}$ $\ast$\\
\hline
      &   &     &       &           &
$(8\aaa_1,8\aaa_1)\subset 8\aaa_2$ & $19$   &
$2_7^{-3},3^{-1}$ $\ast$\\
\hline
      &   &     &       &           &
$(8\aaa_1,\aaa_2)\subset 8\aaa_1\amalg \aaa_2$ & $19$   &
$2_2^{+2},16_3^{-1}$ $\ast$\\
\hline
      &   &     &       &           &
$(\aaa_2,\aaa_2)\subset 2\aaa_2$ & $19$   &
$2_5^{-1},8_{II}^{-2}$ $\ast$\\
\hline
\hline
 $16$ &$10$&$1$ &$D_{10}$&  $16$   & $(\aaa_1,\aaa_1)\subset 2\aaa_1$  & $18$  &
$2_6^{+2},5^{+4}$\\
\hline
      &   &     &       &            &  $(\aaa_1,5\aaa_1)\subset 6\aaa_1$  &$18$  &
$2_6^{+2},5^{-3}$\\
\hline
      &   &     &       &            &  $(5\aaa_1,5\aaa_1)\subset 10\aaa_1$  &$18$  &
$2_6^{+2},5^{+2}$ $\ast$\\
\hline
\end{tabular}
\end{table}

\begin{table}

\begin{tabular}{|c||c|c|c|c|c|c|c|c|}
\hline
 &   &   &  & &  $(6\aaa_1,6\aaa_1)\subset 6\aaa_2$ &  $18$  &
$4_{II}^{-2},3^{+1}$  $\ast$ \\
\hline
 &   &   &  & &  $(6\aaa_1,12\aaa_1)\subset 6\aaa_3$ &  $18$  &
$4_6^{+2}$  $\ast$ \\
\hline
&   &   &  & &  $(\aaa_1,\aaa_1,\aaa_1)\subset 3\aaa_1$ &  $19$  &
$2_{1}^{+1},4_{II}^{+2},3^{+2}$  $\ast$ \\
\hline
&   &   &  & &  $(\aaa_1,\aaa_1,4\aaa_1)\subset 6\aaa_1$ &  $19$  &
$2_{2}^{+2},8_{7}^{+1},3^{+2}$  $\ast$ \\
\hline
&   &   &  & &  $(\aaa_1,\aaa_1,6\aaa_1)\subset 8\aaa_1$ &  $19$  &
$4_{7}^{-3},3^{+1}$  $\ast$ \\
\hline
&   &   &  & &  $(\aaa_1,\aaa_1,12\aaa_1)\subset 14\aaa_1$ &  $19$  &
$2_{2}^{+2},8_{5}^{-1},3^{-1}$  $\ast$ \\
\hline
&   &   &  & &  $(\aaa_1,3\aaa_1,4\aaa_1)\subset 8\aaa_1$ &  $19$  &
$2_{II}^{-2},8_{7}^{+1},3^{-1}$  $\ast$ \\
\hline
&   &   &  & &  $(\aaa_1,3\aaa_1,12\aaa_1)\subset 16\aaa_1$ &  $19$  &
$2_{II}^{-2},8_{5}^{-1}$  $\ast$ \\
\hline

\end{tabular}
\end{table}

\begin{table}

\begin{tabular}{|c||c|c|c|c|c|c|}

\hline
 {\bf n}& $|G|$& $i$&  $G$   &$Deg$& $\rk S$ & $q_S$            \\
\hline
\hline
 $18$&$12$ &$4$ &$D_{12}$&   $(\aaa_1,2\aaa_1)\subset 3\aaa_1$ &  $18$  &
$2_7^{+3},4_7^{+1},3^{+4}$\\
\hline
 &   &   &  &   $(\aaa_1,3\aaa_1)\subset 4\aaa_1$ &  $18$  &
$2_0^{+4},3^{-3}$ $\ast$\\
\hline
&   &    &  &    $(\aaa_1,6\aaa_1)\subset 7\aaa_1$ &  $18$  &
$2_1^{-3},4_7^{+1},3^{+3}$ $\ast$\\
\hline
&   &    &  &    $(2\aaa_1,3\aaa_1)\subset 5\aaa_1$ &  $18$  &
$2_1^{+3},4_7^{+1},3^{-3}$\\
\hline
&   &    &  &    $\left((2\aaa_1,6\aaa_1)\subset 8\aaa_1\right)_I$ &  $18$  &
$4_4^{-2},3^{+3}$ $\ast$ \\
\hline
&   &    &  &    $\left((2\aaa_1,6\aaa_1)\subset 8\aaa_1\right)_{II}$ &  $18$  &
$2_{II}^{-4},3^{+3}$ \\
\hline
&   &    &  &    $(2\aaa_1,6\aaa_1)\subset 2\ddd_4$ &  $18$  &
$2_{II}^{-2},3^{+3}$ $\ast$\\
\hline
\end{tabular}
\end{table}

\begin{table}
\begin{tabular}{|c||c|c|c|c|c|c|c|c|}

\hline
&   &    &  &
$
\left(\begin{array}{rrr}
2\aaa_1 & 5\aaa_1 & (8\aaa_1)_{II}  \\
         & 3\aaa_1 & 3\aaa_3   \\
         &         & 6\aaa_1
\end{array}\right)
\subset 2\aaa_1\amalg 3\aaa_3$
&  $19$  & $2_{II}^{+2},4_1^{+1},3^{+2}$\\
\hline
&   &    &  &
$
\left(\begin{array}{rrr}
2\aaa_1 & 5\aaa_1 & 2\ddd_4  \\
         & 3\aaa_1 & 9\aaa_1   \\
         &         & 6\aaa_1
\end{array}\right)
\subset 3\aaa_1\amalg 2\ddd_4$
&  $19$  & $2_5^{+3},3^{-2}$ $\ast$\\
\hline
&   &    &  &
$
\left(\begin{array}{rrr}
2\aaa_1 & (8\aaa_1)_I & (8\aaa_1)_{II}  \\
         & 6\aaa_1 & 12\aaa_1   \\
         &         & 6\aaa_1
\end{array}\right)
\subset 14\aaa_1$
&  $19$  & $2_{II}^{+2},4_1^{+1},3^{+2}$\\
\hline
&   &    &  &
$
\left(\begin{array}{rrr}
2\aaa_1 & (8\aaa_1)_I & (8\aaa_1)_I  \\
         & 6\aaa_1 & 6\aaa_2   \\
         &         & 6\aaa_1
\end{array}\right)
\subset 2\aaa_1\amalg 6\aaa_2$
&  $19$  & $4_3^{-1},3^{-1},9^{-1}$ $\ast$\\
\hline
&   &    &  &
$
\left(\begin{array}{rrr}
2\aaa_1 & 2\ddd_4 & (8\aaa_1)_I  \\
         & 6\aaa_1 & 12\aaa_1   \\
         &         & 6\aaa_1
\end{array}\right)
\subset 2\ddd_4\amalg 6\aaa_1$
&  $19$  & $4_1^{+1},3^{+2}$ $\ast$\\
\hline
&   &    &  &
$(3\aaa_1,6\aaa_1,6\aaa_1)\subset 15 \aaa_1$
&  $19$  & $2_3^{-1},4_4^{-2},3^{-1}$  $\ast$\\
\hline
&   &    &  &
$(3\aaa_1,6\aaa_1,6\aaa_1)\subset 3\aaa_1\amalg 6\aaa_2 $
&  $19$  & $2_1^{-3},9^{-1}$  $\ast$\\
\hline
&   &    &  &
$
\left(\begin{array}{rrr}
3\aaa_1 & 3\aaa_3 & 9\aaa_1  \\
         & 6\aaa_1 & 12\aaa_1   \\
         &         & 6\aaa_1
\end{array}\right)
\subset 3\aaa_3\amalg 6\aaa_1$
&  $19$  & $2_6^{-2},4_5^{-1},3^{+1}$ $\ast$\\
\hline

\end{tabular}
\end{table}

\begin{table}

\begin{tabular}{|c||c|c|c|c|c|c|c|c|}
\hline
 {\bf n}& $|G|$& $i$&  $G$   & $Deg$& $\rk S$ &$q_S$            \\
\hline
\hline
 $22$& $16$&$11$&$C_2\times D_8$&   $(2\aaa_1,2\aaa_1)\subset 4\aaa_1$
 &  $18$  & $4_6^{+4}$ $\ast$\\
\hline
     &     &    &        &              $(2\aaa_1,4\aaa_1)\subset 6\aaa_1$
     &  $18$  &
$4_7^{+3},8_7^{+1}$ $\ast$ \\
\hline
     &     &    &        &    $(2\aaa_1,8\aaa_1)\subset 10\aaa_1$    &  $18$  &
$4_6^{+4}$\\
\hline
     &     &    &        &    $((4\aaa_1,4\aaa_1)\subset 8\aaa_1)_I$    &  $18$  &
$4_6^{+4}$ $\ast$\\
\hline
     &     &    &        &   $((4\aaa_1,4\aaa_1)\subset 8\aaa_1)_{II}$    &  $18$  &
$2_{II}^{-2},8_6^{-2}$ $\ast$\\
\hline
     &     &    &        &      $(4\aaa_1,8\aaa_1)\subset 12\aaa_1$    &  $18$  &
$2_{II}^{+2},4_7^{+1},8_7^{+1}$ $\ast$\\
\hline
     &     &    &        &           $(4\aaa_1,8\aaa_1)\subset 4\aaa_3$    &  $18$  &
$2_{II}^{-2},4_2^{-2}$ $\ast$\\
\hline

\end{tabular}

\end{table}

\begin{table}

\begin{tabular}{|c||c|c|c|c|c|c|c|c|}
\hline
 {\bf n}& $|G|$& $i$&  $G$ & $\rk S_G$  & $Deg$& $\rk S$ &$q_S$   \\
\hline
\hline
 $34$& $24$&$12$&$\SSS_4$    & $17$ &
$(\aaa_1,\aaa_1)\subset 2\aaa_1$&  $19$  &
$4_1^{+3},3^{+2}$ $\ast$\\
\hline
& & & &  &  $(\aaa_1,2\aaa_1)\subset 3\aaa_1$&  $19$  &
$2_1^{+1},4_{II}^{+2},3^{+2}$ $\ast$\\
\hline
& & & & & $(\aaa_1,3\aaa_1)\subset 4\aaa_1$&  $19$  &
$4_{3}^{-3},3^{-1}$ $\ast$\\
\hline
& & & & & $(\aaa_1,4\aaa_1)\subset 5\aaa_1$&  $19$  &
$2_1^{+1},4_1^{+1},8_7^{+1},3^{+2}$ $\ast$\\
\hline
& & & & & $(\aaa_1,(6\aaa_1)_{I})\subset 7\aaa_1$&  $19$  &
$2_1^{+1},4_{6}^{-2},3^{+1}$ $\ast$\\
\hline
& & & & & $(\aaa_1,8\aaa_1)\subset 9\aaa_1$&  $19$  &
$2_1^{+1},4_{0}^{+2},3^{+2}$ $\ast$\\
\hline
& & & & & $(\aaa_1,12\aaa_1)\subset 13\aaa_1$&  $19$  &
$2_7^{+1},4_{5}^{-1},8_7^{+1},3^{-1}$\\
\hline
\end{tabular}

\end{table}

\vfill

\newpage

\begin{table}
\label{table3}
\caption{Markings $S\subset N_j$ by Niemeier lattices, and lattices $S^\perp_{N_j}$
of degenerations of
codimension $\ge 2$ of K\"ahlerian K3 surfaces with finite
symplectic automorphism groups $G=Clos(G)$
in notations \cite{Nik8}, \cite{Nik9}, \cite{Nik10} for ${\bf n}\ge 12$.}
\end{table}

{\bf n}=12, degeneration $(8\aaa_1,8\aaa_1)\subset 16\aaa_1$:

\begin{center}
\begin{tabular}{|c|c|c|c|c|c|c|}
\hline
$j$                     &        $23$ $\ast$      \\
\hline
$H$                     &         $H_{12,1}$      \\
\hline
orbits of               &  $(\alpha_1,\alpha_4)$  \\
\hline
$(S^\perp_{N_j})^{(2)}$ &$4A_1$        \\
\hline
\end{tabular}
\end{center}

\medskip

{\bf n}=12, degeneration $(8\aaa_1,8\aaa_1)\subset 8\aaa_2$:

\begin{center}
\begin{tabular}{|c|c|c|c|c|c|c|}
\hline
$j$                    &  $18$                          & $22$ $\ast$  \\
\hline
$H$                    &$H_{12,1}$                      & $H_{12,1}$    \\
\hline
orbits of               & $(\alpha_{1,1},\alpha_{2,1})$ &
$(\alpha_{1,1},\alpha_{2,1})$\\
\hline
$(S^\perp_{N_j})^{(2)}$ & $D_4$                         &
$3A_1\oplus A_2$        \\
\hline
\end{tabular}
\end{center}

\medskip

{\bf n}=12, degeneration $(8\aaa_1,\aaa_2)\subset 8\aaa_1\amalg \aaa_2$:

\begin{center}
\begin{tabular}{|c|c|c|c|c|c|c|}
\hline
$j$       & $22$ $\ast$            \\
\hline
$H$       &$H_{12,1}$        \\
\hline
orbits of  &$(\alpha_{1,1},\alpha_{1,5})$,$(\alpha_{1,1},\alpha_{1,8})$,
$(\alpha_{1,1},\alpha_{1,11})$,
$(\alpha_{2,1},\alpha_{1,5})$,$(\alpha_{2,1},\alpha_{1,8})$,
$(\alpha_{2,1},\alpha_{1,11})$             \\
\hline
$(S^\perp_{N_j})^{(2)}$& $2A_1\oplus A_2$        \\
\hline
\end{tabular}
\end{center}

\medskip

{\bf n}=12, degeneration $(\aaa_2,\aaa_2)\subset 2\aaa_2$:

\begin{center}
\begin{tabular}{|c|c|c|c|c|c|c|}
\hline
$j$      & $22$ $\ast$            \\
\hline
$H$      & $H_{12,1}$                 \\
\hline
orbits of& $(\alpha_{1,5},\alpha_{1,8})$,$(\alpha_{1,5},\alpha_{1,11})$,
$(\alpha_{1,8},\alpha_{1,11})$\\
\hline
$(S^\perp_{N_j})^{(2)}$&$A_1\oplus A_2$        \\
\hline
\end{tabular}
\end{center}

\medskip

{\bf n}=16, degeneration $(\aaa_1,\aaa_1)\subset 2\aaa_1$:

\begin{center}
\begin{tabular}{|c|c|c|c|c|c|c|c|c|}
\hline
$j$ & $19$ & $20$ & $20$ & $22$ & $23$         \\
\hline
$H$ & $H_{16,1}$ & $H_{16,1}$ & $H_{16,1}$ & $H_{16,1}$ & $H_{16,1}$ \\
\hline
orbits of &$(\alpha_{1,3},\alpha_{3,3}),...$ &
$(\alpha_{1,1},\alpha_{3,1}),...$ &
$(\alpha_{1,1},\alpha_{4,1})$ & $(\alpha_{1,1},\alpha_{1,2}),... $ &
$(\alpha_1,\alpha_{6}),... $              \\
\hline
$(S^\perp_{N_j})^{(2)}$ & $2A_1$& $\{0\}$ & $\{0\}$ & $\{0\}$ & $2A_1$  \\
\hline
\end{tabular}
\end{center}

\medskip

{\bf n}=16, degeneration $(\aaa_1,5\aaa_1)\subset 6\aaa_1$:

\begin{center}
\begin{tabular}{|c|c|c|c|c|c|c|c|c|}
\hline
$j$ & $19$ & $19$ & $19$ & $19$ & $19$         \\
\hline
$H$ & $H_{16,1}$ & $H_{16,1}$ & $H_{16,1}$ & $H_{16,1}$ & $H_{16,1}$\\
\hline
orbits of &$(\alpha_{1,3},\alpha_{1,1}),...$ &
$(\alpha_{3,3},\alpha_{1,1}),... $ &
$(\alpha_{2,3},\alpha_{1,1}),...$ & $(\alpha_{1,3},\alpha_{2,1}),... $ &
$(\alpha_{2,3},\alpha_{2,1})$              \\
\hline
$(S^\perp_{N_j})^{(2)}$ & $3A_1$& $3A_1$ & $3A_1$ & $3A_1$ & $3A_1$  \\
\hline
\end{tabular}
\end{center}

\begin{center}
\begin{tabular}{|c|c|c|c|c|c|c|c|c|}
\hline
$20$ & $20$ & $20$ & $20$ & $22$         \\
\hline
$H_{16,1}$ & $H_{16,1}$ & $H_{16,1}$ & $H_{16,1}$ & $H_{16,1}$ \\
\hline
$(\alpha_{1,1},\alpha_{1,2}),...$ &
$(\alpha_{2,1},\alpha_{1,2}),...$ &
$(\alpha_{1,1},\alpha_{2,2}),...$ &
$(\alpha_{2,1},\alpha_{2,2}),...$ &
$(\alpha_{1,2},\alpha_{1,3}),...$              \\
\hline
$A_2$& $A_2$ & $A_2$ & $A_2$ & $A_2$  \\
\hline
\end{tabular}
\end{center}

\begin{center}
\begin{tabular}{|c|c|c|c|c|c|c|c|c|}
\hline
$22$ & $22$ & $23$ & $23$          \\
\hline
$H_{16,1}$ & $H_{16,1}$ & $H_{16,1}$ & $H_{16,1}$ \\
\hline
$(\alpha_{2,2},\alpha_{1,3}),...$ & $(\alpha_{1,1},\alpha_{1,3}),...$ &
$(\alpha_{1},\alpha_{2}),...$ & $(\alpha_{19},\alpha_{2}),...$  \\
\hline
$A_2$& $A_2$ & $3A_1$ & $3A_1$ \\
\hline
\end{tabular}
\end{center}

\medskip

{\bf n}=16, degeneration $(5\aaa_1,5\aaa_1)\subset 10\aaa_1$:

\begin{center}
\begin{tabular}{|c|c|c|c|c|c|c|c|c|}
\hline
$j$ & $19$ & $20$ & $20$ & $22$ & $23$ $\ast$        \\
\hline
$H$ & $H_{16,1}$ & $H_{16,1}$ & $H_{16,1}$ & $H_{16,1}$ & $H_{16,1}$     \\
\hline
orbits of &$(\alpha_{1,1},\alpha_{3,1})$ & $(\alpha_{1,2},\alpha_{3,2}),...$ &
$(\alpha_{1,2},\alpha_{4,2})$ & $(\alpha_{1,3},\alpha_{1,5}),...$ &
$(\alpha_{2},\alpha_{3}),...$              \\
\hline
$(S^\perp_{N_j})^{(2)}$ & $D_4$& $A_4$ & $A_4$ & $2A_2$ & $4A_1$  \\
\hline
\end{tabular}
\end{center}

\medskip

{\bf n}=16, degeneration $(5\aaa_1,5\aaa_1)\subset 5\aaa_2$:

\begin{center}
\begin{tabular}{|c|c|c|c|c|c|c|c|c|}
\hline
$j$ & $19$ & $20$ & $20$ & $22$ $\ast$        \\
\hline
$H$ & $H_{16,1}$ & $H_{16,1}$ & $H_{16,1}$ & $H_{16,1}$        \\
\hline
orbits of &$(\alpha_{1,1},\alpha_{2,1}),...$ &
$(\alpha_{1,2},\alpha_{2,2}),...$ &
$(\alpha_{2,2},\alpha_{3,2})$ & $(\alpha_{1,3},\alpha_{2,3}),...$   \\
\hline
$(S^\perp_{N_j})^{(2)}$ & $D_4$& $A_4$ & $A_4$ & $2A_2$  \\
\hline
\end{tabular}
\end{center}

\medskip

{\bf n}=16, degeneration $(\aaa_1,\aaa_1,5\aaa_1)\subset 7\aaa_1$:

\begin{center}
\begin{tabular}{|c|c|c|c|c|c|c|c|c|}
\hline
$j$ & $19$ & $19$ & $19$          \\
\hline
$H$ & $H_{16,1}$ & $H_{16,1}$ & $H_{16,1}$              \\
\hline
orbits of &$(\alpha_{1,3},\alpha_{3,3},\alpha_{1,1}),...$ &
$(\alpha_{3,3},\alpha_{4,3},\alpha_{1,1}),...$ &
$(\alpha_{1,3},\alpha_{3,3},\alpha_{2,1}),...$           \\
\hline
$(S^\perp_{N_j})^{(2)}$ & $2A_1$& $2A_1$ & $2A_1$   \\
\hline
\end{tabular}
\end{center}

\begin{center}
\begin{tabular}{|c|c|c|c|c|c|c|c|c|}
\hline
$20$ & $20$ & $20$ & $20$          \\
\hline
$H_{16,1}$ & $H_{16,1}$ & $H_{16,1}$ & $H_{16,1}$              \\
\hline
$(\alpha_{1,1},\alpha_{3,1},\alpha_{1,2}),...$ &
$(\alpha_{1,1},\alpha_{4,1},\alpha_{1,2}),...$ &
$(\alpha_{1,1},\alpha_{3,1},\alpha_{2,2}),...$ &
$(\alpha_{1,1},\alpha_{4,1},\alpha_{2,2}),...$           \\
\hline
$\{0\}$ & $\{0\}$& $\{0\}$ & $\{0\}$   \\
\hline
\end{tabular}
\end{center}

\begin{center}
\begin{tabular}{|c|c|c|c|c|c|c|c|c|}
\hline
$22$ & $22$ & $23$ & $23$          \\
\hline
$H_{16,1}$ & $H_{16,1}$ & $H_{16,1}$ & $H_{16,1}$              \\
\hline
$(\alpha_{1,1},\alpha_{1,2},\alpha_{1,3}),...$ &
$(\alpha_{1,1},\alpha_{2,2},\alpha_{1,3}),...$ &
$(\alpha_{1},\alpha_{6},\alpha_{2}),...$ &
$(\alpha_{1},\alpha_{19},\alpha_{2}),...$           \\
\hline
$\{0\}$ & $\{0\}$& $2A_1$ & $2A_1$   \\
\hline
\end{tabular}
\end{center}

\medskip

{\bf n}=16, degeneration $(\aaa_1,5\aaa_1,5\aaa_1)\subset 11\aaa_1$:

\begin{center}
\begin{tabular}{|c|c|c|c|c|c|c|c|c|}
\hline
$j$ & $19$ & $19$ & $19$          \\
\hline
$H$ & $H_{16,1}$ & $H_{16,1}$ & $H_{16,1}$              \\
\hline
orbits of &$(\alpha_{1,3},\alpha_{1,1},\alpha_{3,1}),...$ &
$(\alpha_{4,3},\alpha_{1,1},\alpha_{3,1}),...$ &
$(\alpha_{2,3},\alpha_{1,1},\alpha_{3,1}),...$           \\
\hline
$(S^\perp_{N_j})^{(2)}$ & $3A_1$& $3A_1$ & $3A_1$   \\
\hline
\end{tabular}
\end{center}

\begin{center}
\begin{tabular}{|c|c|c|c|c|c|c|c|c|}
\hline
$20$ & $20$ & $20$ & $20$          \\
\hline
$H_{16,1}$ & $H_{16,1}$ & $H_{16,1}$ & $H_{16,1}$       \\
\hline
$(\alpha_{1,1},\alpha_{1,2},\alpha_{3,2}),...$ &
$(\alpha_{2,1},\alpha_{1,2},\alpha_{3,2}),...$ &
$(\alpha_{1,1},\alpha_{1,2},\alpha_{4,2}),...$ &
$(\alpha_{2,1},\alpha_{1,2},\alpha_{4,2}),...$           \\
\hline
$A_2$ & $A_2$& $A_2$ & $A_2$   \\
\hline
\end{tabular}
\end{center}

\begin{center}
\begin{tabular}{|c|c|c|c|c|c|c|c|c|}
\hline
$22$ & $22$ & $23$ & $23$          \\
\hline
$H_{16,1}$ & $H_{16,1}$ & $H_{16,1}$ & $H_{16,1}$              \\
\hline
$(\alpha_{1,2},\alpha_{1,3},\alpha_{1,5}),...$ &
$(\alpha_{2,2},\alpha_{1,3},\alpha_{1,5}),...$ &
$(\alpha_{1},\alpha_{2},\alpha_{3}),...$ &
$(\alpha_{6},\alpha_{2},\alpha_{3}),...$           \\
\hline
$A_2$ & $A_2$ & $3A_1$ & $3A_1$    \\
\hline
\end{tabular}
\end{center}

\medskip

{\bf n}=16, degeneration $(\aaa_1,5\aaa_1,5\aaa_1)\subset
\aaa_1\amalg 5\aaa_2$:

\begin{center}
\begin{tabular}{|c|c|c|c|c|c|c|c|c|}
\hline
$j$ & $19$ & $19$ & $19$          \\
\hline
$H$ & $H_{16,1}$ & $H_{16,1}$ & $H_{16,1}$              \\
\hline
orbits of &$(\alpha_{1,3},\alpha_{1,1},\alpha_{2,1}),...$ &
$(\alpha_{2,3},\alpha_{1,1},\alpha_{2,1}),...$ &
$(\alpha_{3,3},\alpha_{1,1},\alpha_{2,1}),...$           \\
\hline
$(S^\perp_{N_j})^{(2)}$ & $3A_1$& $3A_1$ & $3A_1$   \\
\hline
\end{tabular}
\end{center}

\begin{center}
\begin{tabular}{|c|c|c|c|c|c|c|c|c|}
\hline
$20$ & $20$ & $20$ & $20$          \\
\hline
$H_{16,1}$ & $H_{16,1}$ & $H_{16,1}$ & $H_{16,1}$   \\
\hline
$(\alpha_{1,1},\alpha_{1,2},\alpha_{2,2}),...$ &
$(\alpha_{2,1},\alpha_{1,2},\alpha_{2,2}),...$ &
$(\alpha_{1,1},\alpha_{2,2},\alpha_{3,2}),... $ &
$(\alpha_{2,1},\alpha_{2,2},\alpha_{3,2}),...$           \\
\hline
$A_2$ & $A_2$& $A_2$ & $A_2$   \\
\hline
\end{tabular}
\end{center}

\begin{center}
\begin{tabular}{|c|c|c|c|c|c|c|c|c|}
\hline
$22$ & $22$           \\
\hline
$H_{16,1}$ & $H_{16,1}$               \\
\hline
$(\alpha_{1,1},\alpha_{1,3},\alpha_{2,3}),...$ &
$(\alpha_{1,2},\alpha_{1,3},\alpha_{2,3}),...$   \\
\hline
$A_2$ & $A_2$     \\
\hline
\end{tabular}
\end{center}

\medskip

{\bf n}=16, degeneration $(5\aaa_1,5\aaa_1,5\aaa_1)\subset 15\aaa_1$:

\begin{center}
\begin{tabular}{|c|c|c|c|c|c|c|c|c|}
\hline
$j$ & $19$ & $23$ $\ast$          \\
\hline
$H$ & $H_{16,1}$ & $H_{16,1}$               \\
\hline
orbits of &$(\alpha_{1,1},\alpha_{3,1},\alpha_{4,1}),...$ &
$(\alpha_{2},\alpha_{3},\alpha_{4}),...$            \\
\hline
$(S^\perp_{N_j})^{(2)}$ & $D_4$& $4A_1$    \\
\hline
\end{tabular}
\end{center}

\medskip

{\bf n}=16, degeneration
$(5\aaa_1,5\aaa_1,5\aaa_1)\subset 5\aaa_2\amalg 5\aaa_1$:

\begin{center}
\begin{tabular}{|c|c|c|c|c|c|c|c|c|}
\hline
$j$ & $20$ & $22$ $\ast$          \\
\hline
$H$ & $H_{16,1}$ & $H_{16,1}$               \\
\hline
orbits of &$(\alpha_{1,2},\alpha_{2,2},\alpha_{4,2}),...$ &
$(\alpha_{1,3},\alpha_{2,3},\alpha_{1,5}),...$         \\
\hline
$(S^\perp_{N_j})^{(2)}$ & $A_4$& $2A_2$    \\
\hline
\end{tabular}
\end{center}

\medskip

{\bf n}=16, degeneration
$
\left(\begin{array}{rrr}
5\aaa_1 & 5\aaa_2 & 10\aaa_1  \\
         & 5\aaa_1 & 5\aaa_2   \\
         &         & 5\aaa_1
\end{array}\right)
\subset 5\aaa_3
$:

\begin{center}
\begin{tabular}{|c|c|c|c|c|c|c|c|c|}
\hline
$j$ & $19$ & $20$ $\ast$          \\
\hline
$H$ & $H_{16,1}$ & $H_{16,1}$               \\
\hline
orbits of &$(\alpha_{1,1},\alpha_{2,1},\alpha_{3,1}),...$ &
$(\alpha_{1,2},\alpha_{2,2},\alpha_{3,2}),...$            \\
\hline
$(S^\perp_{N_j})^{(2)}$ & $D_4$& $A_4$    \\
\hline
\end{tabular}
\end{center}

\medskip

{\bf n}=17, degeneration $(\aaa_1,\aaa_1)\subset 2\aaa_1$:

\begin{center}
\begin{tabular}{|c|c|c|c|c|c|c|c|c|}
\hline
$j$ & $19$ & $21$ & $21$ & $23$ $\ast$  & $23$         \\
\hline
$H$ & $H_{17,2}$ & $H_{17,1}$ & $H_{17,2}$ & $H_{17,1}$ & $H_{17,2}$     \\
\hline
orbits of &$(\alpha_{1,4},\alpha_{3,4}),...$ & $(\alpha_{1,1},\alpha_{3,1})$,&
$(\alpha_{1,1},\alpha_{3,1})$ & $(\alpha_5,\alpha_{17}),...$ &
$(\alpha_5,\alpha_{17})$,              \\
 & & $(\alpha_{1,5},\alpha_{3,5})$ & & & $(\alpha_6,\alpha_{10})$ \\
\hline
$(S^\perp_{N_j})^{(2)}$ & $2A_1\oplus A_2$& $A_3$ & $\{0\}$ & $3A_1$ & $2A_1$ \\
\hline
\end{tabular}
\end{center}

\medskip

{\bf n}=17, degeneration $(\aaa_1,3\aaa_1)\subset 4\aaa_1$:

\begin{center}
\begin{tabular}{|c|c|c|c|c|c|c|c|c|}
\hline
$j$ & $23$ $\ast$  & $23$         \\
\hline
$H$ & $H_{17,1}$ & $H_{17,3}$             \\
\hline
orbits of & $(\alpha_{15},\alpha_{1})$ &
$(\alpha_{15},\alpha_{1}),(\alpha_{10},\alpha_6)$   \\
\hline
$(S^\perp_{N_j})^{(2)}$ & $4A_1$ & $A_1$  \\
\hline
\end{tabular}
\end{center}

{\bf n}=17, degeneration $(\aaa_1,4\aaa_1)\subset 5\aaa_1$:

\begin{center}
\begin{tabular}{|c|c|c|c|c|c|c|c|c|}
\hline
$j$ & $19$ & $19$ & $19$ & $19$ & $19$         \\
\hline
$H$ & $H_{17,1}$ & $H_{17,1}$ & $H_{17,2}$ & $H_{17,2}$ & $H_{17,2}$ \\
\hline
orbits of &$(\alpha_{2,3},\alpha_{1,1}),...$ &
$(\alpha_{2,3},\alpha_{2,1})$, &
$(\alpha_{2,3},\alpha_{2,1})$ & $(\alpha_{1,4},\alpha_{2,1}),...$ &
$(\alpha_{2,4},\alpha_{2,1})$              \\
 & & $(\alpha_{2,4},\alpha_{2,1})$ & & &  \\
\hline
$(S^\perp_{N_j})^{(2)}$ & $A_2$& $A_2$ & $D_4$ &
$3A_1\oplus A_2$ & $3A_1\oplus A_2$  \\
\hline
\end{tabular}
\end{center}

\begin{center}
\begin{tabular}{|c|c|c|c|c|c|c|c|c|}
\hline
$21$ & $21$ & $21$ & $21$ & $22$ & $22$         \\
\hline
$H_{17,2}$ & $H_{17,2}$ & $H_{17,1}$ & $H_{17,2}$ & $H_{17,1}$ & $H_{17,1}$  \\
\hline
$(\alpha_{1,1},\alpha_{1,3}),...$ & $(\alpha_{2,1},\alpha_{1,3})$, &
$(\alpha_{1,1},\alpha_{2,3})$, &
$(\alpha_{2,1},\alpha_{2,3})$ & $(\alpha_{1,1},\alpha_{1,3}),...$  &
$(\alpha_{2,1},\alpha_{1,3}),...$ \\
  &$(\alpha_{2,1},\alpha_{3,3})$ & $(\alpha_{3,1},\alpha_{2,3})$ & & &  \\
\hline
$A_1$ & $A_1$ & $A_1$ & $A_1$ & $A_2$ & $A_2$  \\
\hline
\end{tabular}
\end{center}

\begin{center}
\begin{tabular}{|c|c|c|c|c|c|c|c|c|}
\hline
$23$ & $23$ & $23$ $\ast$ & $23$         \\
\hline
$H_{17,1}$ & $H_{17,1}$ & $H_{17,2}$ & $H_{17,3}$             \\
\hline
$(\alpha_{5},\alpha_{2}),...$ &
$(\alpha_{15},\alpha_{2})$ & $(\alpha_{5},\alpha_{2}),...$ &
$(\alpha_{10},\alpha_{2}),...$              \\
\hline
$4A_1$ & $4A_1$ & $3A_1$ & $A_1$  \\
\hline
\end{tabular}
\end{center}

\medskip

{\bf n}=17, degeneration $(\aaa_1,6\aaa_1)\subset 7\aaa_1$:

\begin{center}
\begin{tabular}{|c|c|c|c|c|c|c|c|c|}
\hline
$j$ & $21$ & $21$ & $21$ & $21$ & $22$ $\ast$         \\
\hline
$H$ & $H_{17,1}$ & $H_{17,1}$ & $H_{17,2}$ & $H_{17,2}$ &$H_{17,1}$ \\
\hline
orbits of &$(\alpha_{1,1},\alpha_{2,2}),...$ &
$(\alpha_{2,1},\alpha_{2,2})$, & $(\alpha_{1,1},\alpha_{1,2})$, &
$(\alpha_{2,1},\alpha_{1,2})$ & $(\alpha_{1,1},\alpha_{1,5}),...$  \\
 & & $(\alpha_{2,5},\alpha_{2,2})$ & $(\alpha_{3,1},\alpha_{1,2})$ & &  \\
\hline
$(S^\perp_{N_j})^{(2)}$ & $A_1\oplus A_3$ & $A_1\oplus A_3$ &
$A_1$ & $A_1$ & $A_2$\\
\hline
\end{tabular}
\end{center}

\begin{center}
\begin{tabular}{|c|c|c|c|c|c|c|c|c|}
\hline
$23$ & $23$         \\
\hline
$H_{17,2}$ & $H_{17,2}$             \\
\hline
$(\alpha_{5},\alpha_{1}),...$ &
$(\alpha_{6},\alpha_{1}),...$               \\
\hline
$3A_1$ & $3A_1$  \\
\hline
\end{tabular}
\end{center}

\medskip

{\bf n}=17, degeneration $(\aaa_1,12\aaa_1)\subset 13\aaa_1$:

\begin{center}
\begin{tabular}{|c|c|c|c|c|c|c|c|c|}
\hline
$j$ & $19$ $\ast$ & $19$ & $19$          \\
\hline
$H$ & $H_{17,2}$ & $H_{17,2}$ & $H_{17,2}$              \\
\hline
orbits of &$(\alpha_{2,3},\alpha_{1,1})$ &
$(\alpha_{1,4},\alpha_{1,1}),...$ & $(\alpha_{2,4},\alpha_{1,1})$    \\
\hline
$(S^\perp_{N_j})^{(2)}$ & $D_4$ & $3A_1\oplus A_2$ & $3A_1\oplus A_2$  \\
\hline
\end{tabular}
\end{center}

\begin{center}
\begin{tabular}{|c|c|c|c|c|c|c|c|c|}
\hline
$21$ & $21$ & $23$ & $23$          \\
\hline
$H_{17,1}$ & $H_{17,1}$ & $H_{17,1}$ &$H_{17,1}$             \\
\hline
$(\alpha_{1,1},\alpha_{1,2}),...$ &
$(\alpha_{2,1},\alpha_{1,2}),(\alpha_{2,5},\alpha_{1,2})$ &
$(\alpha_{5},\alpha_{3}),...$ & $(\alpha_{15},\alpha_3)$   \\
\hline
 $A_1\oplus A_3$ & $A_1\oplus A_3$ & $4A_1$ & $4A_1$  \\
\hline
\end{tabular}
\end{center}

\medskip

{\bf n}=17, degeneration $(3\aaa_1,4\aaa_1)\subset 7\aaa_1$:

\begin{center}
\begin{tabular}{|c|c|c|c|c|c|c|c|c|}
\hline
$j$ & $19$ & $19$ & $19$          \\
\hline
$H$ & $H_{17,1}$ & $H_{17,1}$ & $H_{17,2}$              \\
\hline
orbits of &$(\alpha_{1,3},\alpha_{1,1}),...$ &
$(\alpha_{1,3},\alpha_{2,1}),(\alpha_{1,4},\alpha_{2,1})$ &
$(\alpha_{1,3},\alpha_{2,1})$    \\
\hline
$(S^\perp_{N_j})^{(2)}$ & $A_1\oplus A_2$ & $A_1\oplus A_2$ & $A_1\oplus D_4$  \\
\hline
\end{tabular}
\end{center}

\begin{center}
\begin{tabular}{|c|c|c|c|c|c|c|c|c|}
\hline
$21$ & $21$ & $23$ $\ast$ & $23$          \\
\hline
$H_{17,2}$ & $H_{17,2}$ & $H_{17,1}$ &$H_{17,3}$             \\
\hline
$(\alpha_{2,2},\alpha_{1,3}),(\alpha_{2,2},\alpha_{3,3})$ &
$(\alpha_{2,2},\alpha_{2,3})$ &
$(\alpha_{1},\alpha_{2})$& $(\alpha_{1},\alpha_2),...$   \\
\hline
 $A_3$ & $A_3$ & $5A_1$ & $2A_1$  \\
\hline
\end{tabular}
\end{center}

\medskip

{\bf n}=17, degeneration $(3\aaa_1,6\aaa_1)\subset 3\aaa_3$:

\begin{center}
\begin{tabular}{|c|c|c|c|c|c|c|c|c|}
\hline
$j$ & $21$  $\ast$          \\
\hline
$H$ & $H_{17,2}$               \\
\hline
orbits of &$(\alpha_{2,2},\alpha_{1,2})$ \\
\hline
$(S^\perp_{N_j})^{(2)}$ & $A_3$   \\
\hline
\end{tabular}
\end{center}

\medskip

{\bf n}=17, degeneration $(3\aaa_1,12\aaa_1)\subset 15\aaa_1$:

\begin{center}
\begin{tabular}{|c|c|c|c|c|c|c|c|c|}
\hline
$j$ & $19$ & $23$ $\ast$           \\
\hline
$H$ & $H_{17,2}$ & $H_{17,1}$               \\
\hline
orbits of &$(\alpha_{1,3},\alpha_{1,1})$ & $(\alpha_{1},\alpha_{3})$     \\
\hline
$(S^\perp_{N_j})^{(2)}$ & $A_1\oplus D_4$ & $5A_1$  \\
\hline
\end{tabular}
\end{center}

\medskip

{\bf n}=17, degeneration $(4\aaa_1,4\aaa_1)\subset 8\aaa_1$:

\begin{center}
\begin{tabular}{|c|c|c|c|c|c|c|c|c|}
\hline
$j$ & $19$ & $21$ & $23$ $\ast$ & $23$         \\
\hline
$H$ & $H_{17,1}$ & $H_{17,2}$ & $H_{17,2}$ & $H_{17,3}$              \\
\hline
orbits of &$(\alpha_{1,1},\alpha_{3,1}),...$ &
$(\alpha_{1,3},\alpha_{3,3})$ & $(\alpha_{2},\alpha_{4})$ &
$(\alpha_{2},\alpha_3),...$    \\
\hline
$(S^\perp_{N_j})^{(2)}$ & $2A_2$ & $A_3$ & $4A_1$ & $2A_1$  \\
\hline
\end{tabular}
\end{center}

\medskip

{\bf n}=17, degeneration $(4\aaa_1,4\aaa_1)\subset 4\aaa_2$:

\begin{center}
\begin{tabular}{|c|c|c|c|c|c|c|c|c|}
\hline
$j$ & $19$ & $21$ $\ast$ & $22$          \\
\hline
$H$ & $H_{17,1}$ & $H_{17,2}$ & $H_{17,1}$        \\
\hline
orbits of &$(\alpha_{1,1},\alpha_{2,1}),...$ &
$(\alpha_{1,3},\alpha_{2,3}),(\alpha_{2,3},\alpha_{3,3})$ &
$(\alpha_{1,3},\alpha_{2,3})$   \\
\hline
$(S^\perp_{N_j})^{(2)}$ & $2A_2$ & $A_3$ & $2A_2$\\
\hline
\end{tabular}
\end{center}

\medskip

{\bf n}=17, degeneration $(4\aaa_1,6\aaa_1)\subset 10\aaa_1$:

\begin{center}
\begin{tabular}{|c|c|c|c|c|c|c|c|c|}
\hline
$j$ & $21$ & $21$ & $22$ & $23$  $\ast$         \\
\hline
$H$ & $H_{17,2}$ & $H_{17,2}$ & $H_{17,1}$ & $H_{17,2}$   \\
\hline
orbits of &$(\alpha_{1,3},\alpha_{1,2}),(\alpha_{3,3},\alpha_{1,2})$ &
$(\alpha_{2,3},\alpha_{1,2})$ & $(\alpha_{1,3},\alpha_{1,5})$ &
$(\alpha_{2},\alpha_1),...$    \\
\hline
$(S^\perp_{N_j})^{(2)}$ & $A_3$ & $A_3$ & $2A_2$ & $4A_1$  \\
\hline
\end{tabular}
\end{center}

\medskip

{\bf n}=17, degeneration $(4\aaa_1,12\aaa_1)\subset 16\aaa_1$:

\begin{center}
\begin{tabular}{|c|c|c|c|c|c|c|c|c|}
\hline
$j$ & $23$  $\ast$         \\
\hline
$H$ & $H_{17,1}$              \\
\hline
orbits of &  $(\alpha_{2},\alpha_3)$    \\
\hline
$(S^\perp_{N_j})^{(2)}$ & $5A_1$  \\
\hline
\end{tabular}
\end{center}

\medskip

{\bf n}=17, degeneration $(4\aaa_1,12\aaa_1)\subset 4\ddd_4$:

\begin{center}
\begin{tabular}{|c|c|c|c|c|c|c|c|c|}
\hline
$j$ & $19$ $\ast$    \\
\hline
$H$ & $H_{17,2}$    \\
\hline
orbits of &$(\alpha_{2,1},\alpha_{1,1})$   \\
\hline
$(S^\perp_{N_j})^{(2)}$ & $A_2\oplus D_4$ \\
\hline
\end{tabular}
\end{center}

\medskip

{\bf n}=17, degeneration $(6\aaa_1,6\aaa_1)\subset 12\aaa_1$:

\begin{center}
\begin{tabular}{|c|c|c|c|c|c|c|c|c|}
\hline
$j$ & $23$ $\ast$    \\
\hline
$H$ & $H_{17,2}$    \\
\hline
orbits of &$(\alpha_{1},\alpha_{7})$   \\
\hline
$(S^\perp_{N_j})^{(2)}$ & $4A_1$ \\
\hline
\end{tabular}
\end{center}

\medskip

{\bf n}=17, degeneration $(6\aaa_1,6\aaa_1)\subset 6\aaa_2$:

\begin{center}
\begin{tabular}{|c|c|c|c|c|c|c|c|c|}
\hline
$j$ & $22$ $\ast$    \\
\hline
$H$ & $H_{17,1}$    \\
\hline
orbits of &$(\alpha_{1,5},\alpha_{2,5})$   \\
\hline
$(S^\perp_{N_j})^{(2)}$ & $2A_2$ \\
\hline
\end{tabular}
\end{center}

\medskip

{\bf n}=17, degeneration $(6\aaa_1,12\aaa_1)\subset 6\aaa_3$:

\begin{center}
\begin{tabular}{|c|c|c|c|c|c|c|c|c|}
\hline
$j$ & $21$ $\ast$    \\
\hline
$H$ & $H_{17,1}$    \\
\hline
orbits of &$(\alpha_{2,2},\alpha_{1,2})$   \\
\hline
$(S^\perp_{N_j})^{(2)}$ & $2A_3$ \\
\hline
\end{tabular}
\end{center}

\medskip

{\bf n}=17, degeneration $(\aaa_1,\aaa_1,\aaa_1)\subset 3\aaa_1$:

\begin{center}
\begin{tabular}{|c|c|c|c|c|c|c|c|c|}
\hline
$j$ & $19$ & $23$ $\ast$  \\
\hline
$H$ & $H_{17,2}$ & $H_{17,1}$ \\
\hline
orbits of & $(\alpha_{1,4},\alpha_{3,4},\alpha_{4,4})$ &
$(\alpha_{5},\alpha_{17},\alpha_{21}),...$   \\
\hline
$(S^\perp_{N_j})^{(2)}$ & $A_1\oplus A_2$ & $2A_1$   \\
\hline
\end{tabular}
\end{center}

\medskip

{\bf n}=17, degeneration $(\aaa_1,\aaa_1,4\aaa_1)\subset 6\aaa_1$:

\begin{center}
\begin{tabular}{|c|c|c|c|c|c|c|c|c|}
\hline
$j$ & $19$ & $21$ & $21$ \\
\hline
$H$ & $H_{17,2}$ & $H_{17,2}$ & $H_{17,2}$  \\
\hline
orbits of & $(\alpha_{1,4},\alpha_{3,4},\alpha_{2,1}),...$ &
$(\alpha_{1,1},\alpha_{3,1},\alpha_{1,3}),(\alpha_{1,1},\alpha_{3,1},\alpha_{3,3})$ &
$(\alpha_{1,1},\alpha_{3,1},\alpha_{2,3})$    \\
\hline
$(S^\perp_{N_j})^{(2)}$ & $2A_1\oplus A_2$ & $\{0\}$ & $\{0\}$  \\
\hline
\end{tabular}
\end{center}

\begin{center}
\begin{tabular}{|c|c|c|c|c|c|c|c|c|}
\hline
 $23$ $\ast$ & $23$ \\
\hline
$H_{17,1}$ & $H_{17,2}$  \\
\hline
$(\alpha_{5},\alpha_{17},\alpha_{2}),...$ &
$(\alpha_{5},\alpha_{17},\alpha_{2}),...$    \\
\hline
$3A_1$ & $2A_1$  \\
\hline
\end{tabular}
\end{center}

\medskip

{\bf n}=17, degeneration $(\aaa_1,\aaa_1,6\aaa_1)\subset 8\aaa_1$:

\begin{center}
\begin{tabular}{|c|c|c|c|c|c|c|c|c|}
\hline
$j$ & $21$ & $23$ $\ast$ \\
\hline
$H$ & $H_{17,2}$ & $H_{17,2}$ \\
\hline
orbits of & $(\alpha_{1,1},\alpha_{3,1},\alpha_{1,2})$ &
$(\alpha_{6},\alpha_{10},\alpha_{1}),(\alpha_{5},\alpha_{17},\alpha_{7})$\\
\hline
$(S^\perp_{N_j})^{(2)}$ &$\{0\}$ & $2A_1$  \\
\hline
\end{tabular}
\end{center}

\medskip

{\bf n}=17, degeneration $(\aaa_1,\aaa_1,12\aaa_1)\subset 14\aaa_1$:

\begin{center}
\begin{tabular}{|c|c|c|c|c|c|c|c|c|}
\hline
$j$ & $19$ & $21$ & $23$ $\ast$ \\
\hline
$H$ & $H_{17,2}$ & $H_{17,1}$ & $H_{17,1}$ \\
\hline
orbits of & $(\alpha_{1,4},\alpha_{3,4},\alpha_{1,1}),...$ &
$(\alpha_{1,1},\alpha_{3,1},\alpha_{1,2}),(\alpha_{1,5},\alpha_{3,5},\alpha_{1,2})$ &
$(\alpha_{5},\alpha_{17},\alpha_{3}),...$     \\
\hline
$(S^\perp_{N_j})^{(2)}$ &$2A_1\oplus A_2$ & $A_3$ & $3A_1$  \\
\hline
\end{tabular}
\end{center}

\medskip

{\bf n}=17, degeneration $(\aaa_1,3\aaa_1,4\aaa_1)\subset 8\aaa_1$:

\begin{center}
\begin{tabular}{|c|c|c|c|c|c|c|c|c|}
\hline
$j$ & $23$ $\ast$ & $23$ \\
\hline
$H$ & $H_{17,1}$ & $H_{17,3}$ \\
\hline
orbits of & $(\alpha_{15},\alpha_{1},\alpha_{2})$ &
$(\alpha_{15},\alpha_{1},\alpha_{2})$     \\
\hline
$(S^\perp_{N_j})^{(2)}$ & $4A_1$ & $A_1$  \\
\hline
\end{tabular}
\end{center}

\medskip

{\bf n}=17, degeneration $(\aaa_1,3\aaa_1,12\aaa_1)\subset 16\aaa_1$:

\begin{center}
\begin{tabular}{|c|c|c|c|c|c|c|c|c|}
\hline
$j$ & $23$ $\ast$  \\
\hline
$H$ & $H_{17,1}$   \\
\hline
orbits of & $(\alpha_{15},\alpha_{1},\alpha_{3})$      \\
\hline
$(S^\perp_{N_j})^{(2)}$ & $4A_1$   \\
\hline
\end{tabular}
\end{center}

\medskip

\medskip

{\bf n}=17, degeneration $(\aaa_1,4\aaa_1,4\aaa_1)\subset 9\aaa_1$:

\begin{center}
\begin{tabular}{|c|c|c|c|c|c|c|c|c|}
\hline
$j$ & $19$ & $21$ & $21$ \\
\hline
$H$ & $H_{17,1}$ & $H_{17,2}$ & $H_{17,2}$  \\
\hline
orbits of & $(\alpha_{2,3},\alpha_{1,1},\alpha_{3,1}),...$ &
$(\alpha_{1,1},\alpha_{1,3},\alpha_{3,3}),(\alpha_{3,1},\alpha_{1,3},\alpha_{3,3})$ &
$(\alpha_{2,1},\alpha_{1,3},\alpha_{3,3})$    \\
\hline
$(S^\perp_{N_j})^{(2)}$ & $A_2$ & $A_1$ & $A_1$  \\
\hline
\end{tabular}
\end{center}

\begin{center}
\begin{tabular}{|c|c|c|c|c|c|c|c|c|}
\hline
 $23$ $\ast$ & $23$ \\
\hline
$H_{17,2}$ & $H_{17,3}$  \\
\hline
$(\alpha_{5},\alpha_{2},\alpha_{4}),...$ &
$(\alpha_{10},\alpha_{2},\alpha_{3}),...$    \\
\hline
$3A_1$ & $A_1$  \\
\hline
\end{tabular}
\end{center}

\medskip

{\bf n}=17, degeneration $(\aaa_1,4\aaa_1,4\aaa_1)\subset \aaa_1\amalg 4\aaa_2$:

\begin{center}
\begin{tabular}{|c|c|c|c|c|c|c|c|c|}
\hline
$j$ & $19$ & $21$ & $21$ & $22$\ $\ast$ \\
\hline
$H$ & $H_{17,1}$ & $H_{17,2}$ & $H_{17,2}$ & $H_{17,1}$  \\
\hline
orbits of & $(\alpha_{2,3},\alpha_{1,1},\alpha_{2,1}),...$ &
$(\alpha_{1,1},\alpha_{1,3},\alpha_{2,3}),...$ &
$(\alpha_{2,1},\alpha_{1,3},\alpha_{2,3})$, &
$(\alpha_{1,1},\alpha_{1,3},\alpha_{2,3})$  \\
&  &  &$(\alpha_{2,1},\alpha_{2,3},\alpha_{3,3})$ &    \\
\hline
$(S^\perp_{N_j})^{(2)}$ & $A_2$ & $A_1$ & $A_1$ & $A_2$  \\
\hline
\end{tabular}
\end{center}

\medskip

{\bf n}=17, degeneration $(\aaa_1,4\aaa_1,6\aaa_1)\subset 11\aaa_1$:

\begin{center}
\begin{tabular}{|c|c|c|c|c|c|c|c|c|}
\hline
$j$ & $21$ & $21$ & $21$ & $21$ \\
\hline
$H$ & $H_{17,2}$ & $H_{17,2}$ & $H_{17,2}$ & $H_{17,2}$  \\
\hline
orbits of & $(\alpha_{1,1},\alpha_{1,3},\alpha_{1,2}),...$ &
$(\alpha_{2,1},\alpha_{1,3},\alpha_{1,2})$, &
$(\alpha_{1,1},\alpha_{2,3},\alpha_{1,2})$, &
$(\alpha_{2,1},\alpha_{2,3},\alpha_{1,2})$  \\
&  &  $(\alpha_{2,1},\alpha_{3,3},\alpha_{1,2})$
& $(\alpha_{3,1},\alpha_{2,3},\alpha_{1,2})$ &    \\
\hline
$(S^\perp_{N_j})^{(2)}$ & $A_1$ & $A_1$ & $A_1$ & $A_1$  \\
\hline
\end{tabular}
\end{center}

\begin{center}
\begin{tabular}{|c|c|c|c|c|c|c|c|c|}
\hline
$22$ & $22$ & $23$ & $23$ \\
\hline
$H_{17,1}$ & $H_{17,1}$ & $H_{17,2}$ & $H_{17,2}$  \\
\hline
$(\alpha_{1,1},\alpha_{1,3},\alpha_{1,5}),...$ &
$(\alpha_{2,1},\alpha_{1,3},\alpha_{1,5}),...$ &
$(\alpha_{5},\alpha_{2},\alpha_{1}),...$ &
$(\alpha_{6},\alpha_{2},\alpha_{1}),...$  \\
\hline
$A_2$ & $A_2$ & $3A_1$ & $3A_1$  \\
\hline
\end{tabular}
\end{center}

\medskip

{\bf n}=17, degeneration
$(\aaa_1,4\aaa_1,12\aaa_1)\subset \aaa_1\amalg 4\ddd_4$:

\begin{center}
\begin{tabular}{|c|c|c|c|c|c|c|c|c|}
\hline
$j$ & $19$ $\ast$ & $19$ & $19$ \\
\hline
$H$ & $H_{17,2}$ & $H_{17,2}$ & $H_{17,2}$  \\
\hline
orbits of & $(\alpha_{2,3},\alpha_{2,1},\alpha_{1,1})$ &
$(\alpha_{1,4},\alpha_{1,2},\alpha_{1,1}),...$  &
$(\alpha_{2,4},\alpha_{1,2},\alpha_{1,1})$    \\
\hline
$(S^\perp_{N_j})^{(2)}$ & $D_4$ & $A_2\oplus 3A_1$ & $A_2\oplus 3A_1$  \\
\hline
\end{tabular}
\end{center}

{\bf n}=17, degeneration $(\aaa_1,6\aaa_1,6\aaa_1)\subset \aaa_1\amalg 6\aaa_2$:

\begin{center}
\begin{tabular}{|c|c|c|c|c|c|c|c|c|}
\hline
$j$ & $22$\ $\ast$  \\
\hline
$H$ & $H_{17,1}$  \\
\hline
orbits of & $(\alpha_{1,1},\alpha_{1,5},\alpha_{2,5})$   \\
\hline
$(S^\perp_{N_j})^{(2)}$ & $A_2$  \\
\hline
\end{tabular}
\end{center}

\medskip

{\bf n}=17, degeneration $(\aaa_1,6\aaa_1,12\aaa_1)\subset
\aaa_1\amalg 6\aaa_3$:

\begin{center}
\begin{tabular}{|c|c|c|c|c|c|c|c|c|}
\hline
$j$ & $21$ & $21$ \\
\hline
$H$ & $H_{17,1}$ & $H_{17,1}$  \\
\hline
orbits of & $(\alpha_{1,1},\alpha_{2,2},\alpha_{1,2}),...$ &
$(\alpha_{2,1},\alpha_{2,2},\alpha_{1,2}),(\alpha_{2,5},\alpha_{2,2},\alpha_{1,2})$\\
\hline
$(S^\perp_{N_j})^{(2)}$ & $A_1\oplus A_3$ & $A_1\oplus A_3$ \\
\hline
\end{tabular}
\end{center}

\medskip

{\bf n}=17, degeneration $(3\aaa_1,4\aaa_1,4\aaa_1)\subset 11\aaa_1$:

\begin{center}
\begin{tabular}{|c|c|c|c|c|c|c|c|c|}
\hline
$j$ & $19$& $21$ & $23$ $\ast$ \\
\hline
$H$ & $H_{17,1}$ & $H_{17,2}$ & $H_{17,3}$  \\
\hline
orbits of & $(\alpha_{1,3},\alpha_{1,1},\alpha_{3,1}),...$ &
$(\alpha_{2,2},\alpha_{1,3},\alpha_{3,3})$  &
$(\alpha_{1},\alpha_{2},\alpha_{3}),...$    \\
\hline
$(S^\perp_{N_j})^{(2)}$ & $A_1\oplus A_2$ & $A_3$ & $2A_1$  \\
\hline
\end{tabular}
\end{center}

\medskip

{\bf n}=17, degeneration $(3\aaa_1,4\aaa_1,4\aaa_1)\subset 3\aaa_1\amalg 4\aaa_2$:

\begin{center}
\begin{tabular}{|c|c|c|c|c|c|c|c|c|}
\hline
$j$ & $19$ & $21$ $\ast$  \\
\hline
$H$ & $H_{17,1}$ & $H_{17,2}$  \\
\hline
orbits of & $(\alpha_{1,3},\alpha_{1,1},\alpha_{2,1}),...$ &
$(\alpha_{2,2},\alpha_{1,3},\alpha_{2,3}),(\alpha_{2,2},\alpha_{3,3},\alpha_{2,3})$ \\
\hline
$(S^\perp_{N_j})^{(2)}$ & $A_1\oplus A_2$ & $A_3$ \\
\hline
\end{tabular}
\end{center}

\medskip

{\bf n}=17, degeneration $(3\aaa_1,4\aaa_1,12\aaa_1)\subset 3\aaa_1\amalg 4\ddd_4$:

\begin{center}
\begin{tabular}{|c|c|c|c|c|c|c|c|c|}
\hline
$j$ & $19$ $\ast$  \\
\hline
$H$ & $H_{17,2}$   \\
\hline
orbits of & $(\alpha_{1,3},\alpha_{2,1},\alpha_{1,1})$      \\
\hline
$(S^\perp_{N_j})^{(2)}$ & $A_1\oplus D_4$   \\
\hline
\end{tabular}
\end{center}

\medskip

{\bf n}=17, degeneration $(3\aaa_1,6\aaa_1,4\aaa_1)\subset
3\aaa_3\amalg 4\aaa_1$:

\begin{center}
\begin{tabular}{|c|c|c|c|c|c|c|c|c|}
\hline
$j$ & $21$ & $21$ \\
\hline
$H$ & $H_{17,2}$ & $H_{17,2}$  \\
\hline
orbits of & $(\alpha_{2,2},\alpha_{1,2},\alpha_{1,3}),
(\alpha_{2,2},\alpha_{1,2},\alpha_{3,3})$ &
$(\alpha_{2,2},\alpha_{1,2},\alpha_{2,3})$ \\
\hline
$(S^\perp_{N_j})^{(2)}$ & $A_3$ & $A_3$ \\
\hline
\end{tabular}
\end{center}

\medskip

{\bf n}=17, degeneration $(4\aaa_1,4\aaa_1,4\aaa_1)\subset 12\aaa_1$:

\begin{center}
\begin{tabular}{|c|c|c|c|c|c|c|c|c|}
\hline
$j$ & $19$ & $23$ $\ast$ \\
\hline
$H$ & $H_{17,1}$ & $H_{17,3}$  \\
\hline
orbits of & $(\alpha_{1,1},\alpha_{3,1},\alpha_{4,1})$ &
$(\alpha_{2},\alpha_{3},\alpha_{8})$ \\
\hline
$(S^\perp_{N_j})^{(2)}$ & $2A_2$ & $2A_1$ \\
\hline
\end{tabular}
\end{center}

\medskip

{\bf n}=17, degeneration
$
\left(\begin{array}{rrr}
4\aaa_1 & 4\aaa_2 & 8\aaa_1  \\
         & 4\aaa_1 & 4\aaa_2   \\
         &         & 4\aaa_1
\end{array}\right)
\subset 4\aaa_3
$:

\begin{center}
\begin{tabular}{|c|c|c|c|c|c|c|c|c|}
\hline
$j$ & $19$ & $21$ $\ast$ \\
\hline
$H$ & $H_{17,1}$ & $H_{17,2}$  \\
\hline
orbits of & $(\alpha_{1,1},\alpha_{2,1},\alpha_{3,1})$ &
$(\alpha_{1,3},\alpha_{2,3},\alpha_{3,3})$ \\
\hline
$(S^\perp_{N_j})^{(2)}$ & $2A_2$ & $A_3$ \\
\hline
\end{tabular}
\end{center}

\medskip

{\bf n}=17, degeneration $(4\aaa_1,4\aaa_1,6\aaa_1)\subset 14\aaa_1$:

\begin{center}
\begin{tabular}{|c|c|c|c|c|c|c|c|c|}
\hline
$j$ & $21$ & $23$ $\ast$ \\
\hline
$H$ & $H_{17,2}$ & $H_{17,2}$  \\
\hline
orbits of & $(\alpha_{1,3},\alpha_{3,3},\alpha_{1,2})$ &
$(\alpha_{2},\alpha_{4},\alpha_{1}),(\alpha_2,\alpha_4,\alpha_7)$ \\
\hline
$(S^\perp_{N_j})^{(2)}$ & $A_3$ & $4A_1$ \\
\hline
\end{tabular}
\end{center}

\medskip

{\bf n}=17, degeneration $(4\aaa_1,4\aaa_1,6\aaa_1)\subset
4\aaa_2\amalg 6\aaa_1$:

\begin{center}
\begin{tabular}{|c|c|c|c|c|c|c|c|c|}
\hline
$j$ & $21$ & $22$ $\ast$ \\
\hline
$H$ & $H_{17,2}$ & $H_{17,1}$  \\
\hline
orbits of & $(\alpha_{1,3},\alpha_{2,3},\alpha_{1,2}),
(\alpha_{2,3},\alpha_{3,3},\alpha_{1,2})$ &
$(\alpha_{1,3},\alpha_{2,3},\alpha_{1,5}),
(\alpha_{1,3},\alpha_{2,3},\alpha_{2,5})$ \\
\hline
$(S^\perp_{N_j})^{(2)}$ & $A_3$ & $2A_2$ \\
\hline
\end{tabular}
\end{center}

\medskip

{\bf n}=17, degeneration $(4\aaa_1,6\aaa_1,6\aaa_1)\subset 16\aaa_1$:

\begin{center}
\begin{tabular}{|c|c|c|c|c|c|c|c|c|}
\hline
$j$ & $23$ $\ast$ \\
\hline
$H$ & $H_{17,2}$  \\
\hline
orbits of & $(\alpha_{2},\alpha_{1},\alpha_{7}),(\alpha_4,\alpha_1,\alpha_7)$\\
\hline
$(S^\perp_{N_j})^{(2)}$ &  $4A_1$ \\
\hline
\end{tabular}
\end{center}

\medskip

{\bf n}=17, degeneration $(4\aaa_1,6\aaa_1,6\aaa_1)\subset 4\aaa_1\amalg 6\aaa_2$:

\begin{center}
\begin{tabular}{|c|c|c|c|c|c|c|c|c|}
\hline
$j$ & $22$ $\ast$ \\
\hline
$H$ & $H_{17,1}$  \\
\hline
orbits of & $(\alpha_{1,3},\alpha_{1,5},\alpha_{2,5}),
(\alpha_{2,3},\alpha_{1,5},\alpha_{2,5})$ \\
\hline
$(S^\perp_{N_j})^{(2)}$ &  $2A_2$ \\
\hline
\end{tabular}
\end{center}

%%%%%%%%%%%%%%%%%%%%%%%%%%%%%%%%%%%%%%%%%%
%%%%%%%%%%%%%%%%%%%%%%%%%%%%%%%%%%%%%%%%%%
%%%%%%%%%%%%%%%%%%%%%%%%%%%%%%%%%%%%%%%%%%

%%%%%%%%%%%%%%%%%%%%%%%%%%%%%%%%%%%
%%%%%%%%%%%%%%%%%%%%%%%%%%%%%%%%%%%%%%
%%%%%%%%%%%%%%%%%%%%%%%%%%%%%%%%%%%%
%%%%%%%%%%%%%%%%%%%%%%%%%%%%%%%%%%%%%%

\medskip

{\bf n}=18, degeneration $(\aaa_1,2\aaa_1)\subset 3\aaa_1$:

\begin{center}
\begin{tabular}{|c|c|c|c|c|c|c|c|c|c|}
\hline
$j$       & $12$   & $12$ & $12$ & $18$ & $18$ & $18$     \\
\hline
$H$       & $H_{18,1}$ & $H_{18,1}$ &  $H_{18,1}$&
$H_{18,1}$ & $H_{18,1}$ & $H_{18,1}$               \\
\hline
orbits of & $(\alpha_{2,1},\alpha_{1,1})$ &
$(\alpha_{4,1},\alpha_{1,1})$ &$(\alpha_{2,1},\alpha_{3,1})$ &
$(\alpha_{3,1},\alpha_{1,1})$ & $(\alpha_{2,5},\alpha_{1,1})$ &
$(\alpha_{2,5},\alpha_{2,1})$ \\
\hline
$(S^\perp_{N_j})^{(2)}$& $A_1$ & $ A_1$& $A_1$& $A_2$ & $A_1$ & $A_1$  \\
\hline
\end{tabular}
\end{center}

\begin{center}
\begin{tabular}{|c|c|c|c|c|c|c|c|c|c|}
\hline
$19$       & $19$   & $19$ & $19$ & $21$            \\
\hline
$H_{18,1}$ & $H_{18,2}$ &  $H_{18,2}$& $H_{18,2}$ & $H_{18,1}$\\
\hline
$(\alpha_{2,3},\alpha_{2,1}),...$ &
$(\alpha_{4,3},\alpha_{1,3})$ &$(\alpha_{2,3},\alpha_{2,4})$ &
$(\alpha_{4,3},\alpha_{2,4})$ &
$(\alpha_{2,2},\alpha_{1,1}), (\alpha_{2,1},\alpha_{1,2})$\\
\hline
 $A_2$ & $ A_1$& $A_1$& $A_1$ & $A_1$  \\
\hline
\end{tabular}
\end{center}

\begin{center}
\begin{tabular}{|c|c|c|c|c|c|c|c|c|c|}
\hline
$22$       & $22$   & $23$             \\
\hline
$H_{18,1}$ & $H_{18,1}$ &  $H_{18,1}$                \\
\hline
$(\alpha_{1,1},\alpha_{1,7}),(\alpha_{2,1},\alpha_{2,7})$ &
$(\alpha_{2,1},\alpha_{1,7}),(\alpha_{1,1},\alpha_{2,7})$ &
$(\alpha_1,\alpha_{15}),(\alpha_6,\alpha_{15}),
(\alpha_{1},\alpha_{16}),(\alpha_6,\alpha_{16})$ \\
\hline
 $\{0\}$ & $\{0\}$& $A_1$   \\
\hline
\end{tabular}
\end{center}

\medskip

{\bf n}=18, degeneration $(\aaa_1,3\aaa_1)\subset 4\aaa_1$:

\begin{center}
\begin{tabular}{|c|c|c|c|c|c|c|c|c|c|}
\hline
$j$       & $18$   & $21$ \ $\ast$            \\
\hline
$H$       & $H_{18,1}$ & $H_{18,1}$          \\
\hline
orbits of & $(\alpha_{3,1},\alpha_{3,2})$ &
$(\alpha_{2,2},\alpha_{2,3}),(\alpha_{2,1},\alpha_{2,5})$\\
\hline
$(S^\perp_{N_j})^{(2)}$& $2A_1\oplus A_2$ & $3A_1$ \\
\hline
\end{tabular}
\end{center}

\medskip

{\bf n}=18, degeneration $(\aaa_1,6\aaa_1)\subset 7\aaa_1$:

\begin{center}
\begin{tabular}{|c|c|c|c|c|c|c|c|c|c|}
\hline
$j$       & $12$   & $12$ & $12$ & $12$ & $18$ & $18$      \\
\hline
$H$       & $H_{18,1}$ & $H_{18,1}$ & $H_{18,1}$ & $H_{18,1}$ &
$H_{18,1}$ & $H_{18,1}$ \\
\hline
orbits of & $(\alpha_{2,1},\alpha_{1,2})$ & $(\alpha_{4,1},\alpha_{1,2})$ &
$(\alpha_{2,1},\alpha_{3,2})$ & $(\alpha_{4,1},\alpha_{3,2})$ &
$(\alpha_{3,1},\alpha_{1,2})$ & $(\alpha_{2,5},\alpha_{1,2})$ \\
\hline
$(S^\perp_{N_j})^{(2)}$& $3A_1$ & $3A_1$ & $3A_1$ & $3A_1$ &
$2A_1\oplus A_2$ & $3A_1$ \\
\hline
\end{tabular}
\end{center}

\begin{center}
\begin{tabular}{|c|c|c|c|c|c|c|c|c|c|}
\hline
$18$   & $18$           &   $19$\ $\ast$ & $19$ & $19$ & $19$        \\
\hline
$H_{18,1}$ & $H_{18,1}$ & $H_{18,1}$ & $H_{18,2}$ &
$H_{18,2}$ & $H_{18,2}$          \\
\hline
$(\alpha_{3,1},\alpha_{2,2})$ & $(\alpha_{2,5},\alpha_{2,2})$ &
$(\alpha_{2,3},\alpha_{1,1}),(\alpha_{2,4},\alpha_{1,1})$ &
$(\alpha_{2,3},\alpha_{1,1})$ &
$(\alpha_{4,3},\alpha_{1,1})$ & $(\alpha_{2,3},\alpha_{1,4})$\\
\hline
$2A_1\oplus A_2$ & $3A_1$ &
$A_2$ & $A_1$ & $A_1$ & $A_1$  \\
\hline
\end{tabular}
\end{center}

\begin{center}
\begin{tabular}{|c|c|c|c|c|c|c|c|c|c|}
\hline
$19$   & $21$ & $21$ & $22$              \\
\hline
$H_{18,2}$ & $H_{18,1}$ & $H_{18,1}$ & $H_{18,1}$           \\
\hline
$(\alpha_{4,3},\alpha_{1,4})$ &
$(\alpha_{2,1},\alpha_{1,3}),(\alpha_{2,2},\alpha_{1,5})$ &
$(\alpha_{2,2},\alpha_{1,3}),(\alpha_{2,1},\alpha_{1,5})$ &
$(\alpha_{1,1},\alpha_{1,5}),(\alpha_{2,1},\alpha_{2,5})$ \\
\hline
$A_1$ & $3A_1$ & $3A_1$ & $\{0\}$  \\
\hline
\end{tabular}
\end{center}

\begin{center}
\begin{tabular}{|c|c|c|c|c|c|c|c|c|c|}
\hline
$22$   & $23$             \\
\hline
$H_{18,1}$ & $H_{18,1}$          \\
\hline
$(\alpha_{2,1},\alpha_{1,5}),(\alpha_{1,1},\alpha_{2,5})$ &
$(\alpha_1,\alpha_2),(\alpha_6,\alpha_2),(\alpha_1,\alpha_8),
(\alpha_6,\alpha_8)$\\
\hline
$\{0\}$  & $A_1$\\
\hline
\end{tabular}
\end{center}

\medskip

{\bf n}=18, degeneration $(2\aaa_1,3\aaa_1)\subset 5\aaa_1$:

\begin{center}
\begin{tabular}{|c|c|c|c|c|c|c|c|c|c|}
\hline
$j$       & $12$   & $12$ & $12$ & $12$ & $18$ & $18$         \\
\hline
$H$       & $H_{18,1}$ & $H_{18,1}$ & $H_{18,1}$ & $H_{18,1}$ &
$H_{18,1}$ & $H_{18,1}$         \\
\hline
orbits of & $(\alpha_{1,1},\alpha_{2,2})$ & $(\alpha_{1,1},\alpha_{4,2})$ &
$(\alpha_{3,1},\alpha_{2,2})$ & $(\alpha_{3,1},\alpha_{4,2})$ &
$(\alpha_{1,1},\alpha_{3,2})$ & $(\alpha_{1,1},\alpha_{1,5})$ \\
\hline
$(S^\perp_{N_j})^{(2)}$ & $A_3$ & $A_3$ & $A_3$ & $A_3$ &
$A_1\oplus A_2$ & $2A_1$  \\
\hline
\end{tabular}
\end{center}

\begin{center}
\begin{tabular}{|c|c|c|c|c|c|c|c|c|c|}
\hline
$18$  & $18$ & $19$ & $19$ & $19$ & $19$              \\
\hline
$H_{18,1}$ & $H_{18,1}$ & $H_{18,1}$ &
$H_{18,2}$ & $H_{18,2}$ & $H_{18,2}$           \\
\hline
$(\alpha_{2,1},\alpha_{3,2})$ & $(\alpha_{2,1},\alpha_{1,5})$ &
$(\alpha_{2,1},\alpha_{1,3}),...$ &
$(\alpha_{1,3},\alpha_{2,1})$ & $(\alpha_{2,4},\alpha_{2,1})$ &
$(\alpha_{1,3},\alpha_{3,1})$  \\
\hline
$A_1\oplus A_2$ & $2A_1$ &  $A_1\oplus A_2$ & $2A_1$ & $A_3$ & $2A_1$\\
\hline
\end{tabular}
\end{center}

\begin{center}
\begin{tabular}{|c|c|c|c|c|c|c|c|c|c|}
\hline
$19$  & $21$   & $21$ & $22$               \\
\hline
$H_{18,2}$ & $H_{18,1}$ & $H_{18,1}$ & $H_{18,1}$   \\
\hline
$(\alpha_{2,4},\alpha_{3,1})$ &
$(\alpha_{1,1},\alpha_{2,3}),(\alpha_{1,2},\alpha_{2,5})$ &
$(\alpha_{1,1},\alpha_{2,5}),(\alpha_{1,2},\alpha_{2,3})$ &
$(\alpha_{1,7},\alpha_{1,2}),(\alpha_{2,7},\alpha_{2,2})$  \\
\hline
$A_3$& $2A_1$ & $2A_1$ & $A_2$  \\
\hline
\end{tabular}
\end{center}

\begin{center}
\begin{tabular}{|c|c|c|c|c|c|c|c|c|c|}
\hline
$22$       & $23$             \\
\hline
$H_{18,1}$ & $H_{18,1}$          \\
\hline
$(\alpha_{2,7},\alpha_{1,2}),(\alpha_{1,7},\alpha_{2,2})$ &
$(\alpha_{15},\alpha_{4}),(\alpha_{16},\alpha_4),
(\alpha_{15},\alpha_5),(\alpha_{16},\alpha_5)$\\
\hline
$A_2$ & $2A_1$ \\
\hline
\end{tabular}
\end{center}

\medskip

{\bf n}=18, degeneration $\left((2\aaa_1,6\aaa_1)\subset
8\aaa_1\right)_I$:

\begin{center}
\begin{tabular}{|c|c|c|c|c|c|c|c|c|c|}
\hline
$j$       & $12$   & $12$ & $12$ & $12$ & $18$ & $18$             \\
\hline
$H$  & $H_{18,1}$ & $H_{18,1}$ & $H_{18,1}$ & $H_{18,1}$ &
$H_{18,1}$ & $H_{18,1}$    \\
\hline
orbits of & $(\alpha_{1,1},\alpha_{1,2})$ & $(\alpha_{1,1},\alpha_{3,2})$ &
$(\alpha_{3,1},\alpha_{1,2})$ & $(\alpha_{3,1},\alpha_{3,2})$ &
$(\alpha_{1,1},\alpha_{1,2})$ & $(\alpha_{1,1},\alpha_{2,2})$ \\
\hline
$(S^\perp_{N_j})^{(2)}$& $A_3$ & $A_3$ & $A_3$ & $A_3$ &
$A_1\oplus A_2$ & $A_1\oplus A_2$ \\
\hline
\end{tabular}
\end{center}

\begin{center}
\begin{tabular}{|c|c|c|c|c|c|c|c|c|c|}
\hline
$18$ & $18$   & $19$\   $\ast$& $19$ & $19$             \\
\hline
$H_{18,1}$ & $H_{18,1}$ & $H_{18,1}$ &$H_{18,2}$ & $H_{18,2}$ \\
\hline
$(\alpha_{2,1},\alpha_{1,2})$ & $(\alpha_{2,1},\alpha_{2,2})$ &
$(\alpha_{2,5},\alpha_{1,1}),(\alpha_{2,1},\alpha_{1,5})$ &
$(\alpha_{2,4},\alpha_{1,1})$ & $(\alpha_{1,3},\alpha_{1,4})$  \\
\hline
$A_1\oplus A_2$ & $A_1\oplus A_2$ & $2A_2$ & $A_3$ & $2A_1$  \\
\hline
\end{tabular}
\end{center}

\begin{center}
\begin{tabular}{|c|c|c|c|c|c|c|c|c|c|}
\hline
$21$ & $22$   & $22$ & $23$             \\
\hline
$H_{18,1}$ & $H_{18,1}$ & $H_{18,1}$ & $H_{18,1}$ \\
\hline
$(\alpha_{1,1},\alpha_{1,3}),(\alpha_{1,2},\alpha_{1,5})$ &
$(\alpha_{1,7},\alpha_{1,5}),(\alpha_{2,7},\alpha_{2,5})$ &
$(\alpha_{2,7},\alpha_{1,5}),(\alpha_{1,7},\alpha_{2,5})$ &
$(\alpha_{16},\alpha_2)$,\\
&&& $(\alpha_{15},\alpha_8)$ \\
\hline
$2A_1$& $A_2$ & $A_2$ & $2A_1$  \\
\hline
\end{tabular}
\end{center}

\medskip

{\bf n}=18, degeneration $\left((2\aaa_1,6\aaa_1)\subset
8\aaa_1\right)_{II}$:

\begin{center}
\begin{tabular}{|c|c|c|c|c|c|c|c|c|c|}
\hline
$j$       & $19$   & $21$ & $23$              \\
\hline
$H$       & $H_{18,2}$ & $H_{18,1}$ & $H_{18,1}$         \\
\hline
orbits of & $(\alpha_{1,3},\alpha_{1,1})$ &
$(\alpha_{1,1},\alpha_{1,5}),(\alpha_{1,2},\alpha_{1,3})$ &
$(\alpha_{15},\alpha_2),(\alpha_{16},\alpha_8)$ \\
\hline
$(S^\perp_{N_j})^{(2)}$& $2A_1$ & $2A_1$ & $2A_1$ \\
\hline
\end{tabular}
\end{center}

\medskip

{\bf n}=18, degeneration $(2\aaa_1,6\aaa_1)\subset 2\ddd_4$:

\begin{center}
\begin{tabular}{|c|c|c|c|c|c|c|c|c|c|}
\hline
$j$       & $19$ \ $\ast$  & $19$             \\
\hline
$H$       & $H_{18,1}$ & $H_{18,2}$          \\
\hline
orbits of & $(\alpha_{2,1},\alpha_{1,1}),(\alpha_{2,5},\alpha_{1,5})$ &
$(\alpha_{2,4},\alpha_{1,4})$\\
\hline
$(S^\perp_{N_j})^{(2)}$& $2A_2$ & $A_3$ \\
\hline
\end{tabular}
\end{center}

\medskip

{\bf n}=18, degeneration $(3\aaa_1,6\aaa_1)\subset 9\aaa_1$:

\begin{center}
\begin{tabular}{|c|c|c|c|c|c|c|c|c|c|}
\hline
$j$       & $12$   & $12$ & $12$ & $18$ & $18$ & $18$             \\
\hline
$H$       & $H_{18,1}$ & $H_{18,1}$ & $H_{18,1}$ & $H_{18,1}$ &
$H_{18,1}$ & $H_{18,1}$         \\
\hline
orbits of & $(\alpha_{2,2},\alpha_{1,2})$ & $(\alpha_{4,2},\alpha_{1,2})$ &
$(\alpha_{2,2},\alpha_{3,2})$ & $(\alpha_{3,2},\alpha_{1,2})$
& $(\alpha_{1,5},\alpha_{1,2})$ & $(\alpha_{1,5},\alpha_{2,2})$ \\
\hline
$(S^\perp_{N_j})^{(2)}$& $D_4$ & $D_4$ & $D_4$ &
$3A_1\oplus A_2$ & $4A_1$ & $4A_1$  \\
\hline
\end{tabular}%
\end{center}

\newpage

%%%%%%%%%%%%%%%%%%%%%%%%%%%%%%%%%%%%%%%%%%
%%%%%%%%%%%%%%%%%%%%%%%%%%%%%%%%%%%%%%%%%%%%
%%%%%%%%%%%%%%%%%%%%%%%%%%%%%%%%%%%%%%%%%%%

\section{Tables of types of degenerations of K\"ahlerian K3
surfaces with symplectic automorphism group $D_8$}
\label{sec3.tables}

Unfortunately, for remaining groups of small order (when ${\bf n}<12$), there are many types
of degenerations.
From our point of view, it is reasonable to consider
classification of degenerations for each of these groups separately.
The number of cases depends on a group very much.
Here, we consider the case of $D_8$.

\begin{table}
\label{table4}
\caption{Types and lattices $S$ of degenerations of
codimension $\ge 2$ of K\"ahlerian K3 surfaces
with symplectic automorphism group $D_8$.}


\end{table}

\newpage

\begin{table}
\label{table5}
\caption{Markings $S\subset N_j$ by Niemeier lattices, and lattices $S^\perp_{N_j}$
of degenerations of
codimension $\ge 2$ of K\"ahlerian K3 surfaces with
symplectic automorphism groups $D_8$
in notations \cite{Nik8}, \cite{Nik9}, \cite{Nik10}.}
\end{table}

\medskip

{\bf n}=10, degeneration $(\aaa_1,\aaa_1)\subset 2\aaa_1$:

\begin{center}
%
\end{center}

\newpage

%%%%%%%%%%%%%%%%%%%%%%%%%%%%%%%%%%%%%%%%%%
%%%%%%%%%%%%%%%%%%%%%%%%%%%%%%%%%%%%%%%%%%%%
%%%%%%%%%%%%%%%%%%%%%%%%%%%%%%%%%%%%%%%%%%%

\section{Tables of types of degenerations of K\"ahlerian K3
surfaces with symplectic automorphism group $(C_2)^3$}
\label{sec4.tables}

Here, we consider the case of $(C_2)^3$. This case contains much less cases, 
and it is much easier than the previous case of $D_8$.

\begin{table}
\label{table6}
\caption{Types and lattices $S$ of degenerations of
codimension $\ge 2$ of K\"ahlerian K3 surfaces
with symplectic automorphism group $(C_2)^3$.}

\begin{tabular}{|c||c|c|c|c|c|c|c|c|}
\hline
 {\bf n}& $|G|$ &  $G$   &$Deg$& $\rk S$ &$q_S$            \\
\hline
\hline
  $9$&$8$& $(C_2)^3$ &
  $((2\aaa_1,2\aaa_1)\subset 4\aaa_1)_I$ & $16$  &
$2_{II}^{+2},4_0^{+4}$ $\ast$\\
\hline
           &       &           &
$((2\aaa_1,2\aaa_1)\subset 4\aaa_1)_{II}$ & $16$   &
$2_{II}^{+4},4_{II}^{+2}$ $\ast$ \\
\hline
    &       &           &
$(2\aaa_1,4\aaa_1)\subset 6\aaa_1$ & $16$   &
$2_{II}^{+4},4_7^{+1},8_1^{+1}$ $\ast$\\
\hline
    &       &           &
$(2\aaa_1,4\aaa_1)\subset 2\aaa_3$ & $16$   &
$2_{II}^{+4},4_{II}^{+2}$ $\ast$\\
\hline
    &       &           &
$(2\aaa_1,8\aaa_1)\subset 10\aaa_1$ & $16$   &
$2_{II}^{-4},4_4^{-2}$ $\ast$\\
\hline
    &       &           &
$(4\aaa_1,4\aaa_1)\subset 8\aaa_1$ & $16$   &
$2_{II}^{+4},4_{II}^{+2}$ \\
\hline
   &       &           &
$(4\aaa_1,8\aaa_1)\subset 4\aaa_3$ & $16$   &
$2_{II}^{+6}$ $\ast$\\
\hline
   &       &           &
$(8\aaa_1,8\aaa_1)\subset 16\aaa_1$ & $16$   &
$2_{II}^{+6}$ $\ast$\\
\hline
     &       &           &
$
\left(\begin{array}{ccc}
 2\aaa_1 & (4\aaa_1)_I & (4\aaa_1)_I  \\
         & 2\aaa_1 & (4\aaa_1)_I   \\
         &         & 2\aaa_1
\end{array}\right)
\subset 6\aaa_1
$ & $17$   & $4_7^{+5}$ $\ast$ \\
\hline
     &       &           &
$((2\aaa_1,2\aaa_1)_I,4\aaa_1)\subset 8\aaa_1$ & $17$   &
$2_{II}^{+2},4_6^{+2},8_1^{+1}$ $\ast$\\
\hline
     &       &           &
$((2\aaa_1,2\aaa_1)_{II},4\aaa_1)\subset 8\aaa_1$ & $17$   &
$2_{II}^{+4},8_7^{+1}$ $\ast$\\
\hline
     &       &           &
$
\left(\begin{array}{ccc}
 2\aaa_1 & (4\aaa_1)_I & 6\aaa_1  \\
         & 2\aaa_1 & 2\aaa_3   \\
         &         & 4\aaa_1
\end{array}\right)
\subset 2\aaa_1\amalg 2\aaa_3
$ & $17$   & $2_{II}^{+2},4_7^{+3}$ $\ast$ \\
\hline
     &       &           &
$((2\aaa_1,2\aaa_1)_I,8\aaa_1)\subset 12\aaa_1$ & $17$   &
$2_{II}^{+2},4_7^{+3}$ $\ast$\\
\hline
     &       &           &
$(2\aaa_1,4\aaa_1,4\aaa_1)\subset 10\aaa_1$ & $17$   &
$2_{II}^{+2},4_7^{+3}$ $\ast$ \\
\hline
     &       &           &
$(2\aaa_1,4\aaa_1,4\aaa_1)\subset 2\aaa_3\amalg 4\aaa_1$ & $17$   &
$2_{II}^{+4},8_7^{+1}$ $\ast$ \\
\hline
     &       &           &
$(2\aaa_1,4\aaa_1,8\aaa_1)\subset 2\aaa_1\amalg 4\aaa_3$ & $17$   &
$2_{II}^{+4},4_7^{+1}$ $\ast$ \\
\hline

     &       &           &
$(4\aaa_1,4\aaa_1,4\aaa_1)\subset 12\aaa_1$ & $17$   &
$2_{II}^{+4},8_7^{+1}$ \\
\hline
     &       &           &
$
\left(\begin{array}{cccc}
 2\aaa_1 & (4\aaa_1)_I & (4\aaa_1)_I  & (4\aaa_1)_I \\
         & 2\aaa_1 & (4\aaa_1)_I & (4\aaa_1)_I  \\
         &         & 2\aaa_1 &    (4\aaa_1)_I \\
         &         &         & 2\aaa_1
\end{array}\right)
\subset 8\aaa_1
$ & $18$   & $4_6^{+4}$ $\ast$ \\
\hline
     &       &           &
$
\left(\begin{array}{cccc}
 2\aaa_1 & (4\aaa_1)_I & (4\aaa_1)_I  & 6\aaa_1 \\
         & 2\aaa_1     & (4\aaa_1)_I  & 6\aaa_1 \\
         &             & 2\aaa_1      & 6\aaa_1 \\
         &             &              & 4\aaa_1
\end{array}\right)
\subset 10\aaa_1
$ & $18$   & $4_5^{+3},8_1^{+1}$ $\ast$ \\
\hline
     &       &     &
$
\left(\begin{array}{cccc}
 2\aaa_1 & (4\aaa_1)_I & (4\aaa_1)_I  & 6\aaa_1 \\
         & 2\aaa_1     & (4\aaa_1)_I  & 6\aaa_1 \\
         &             & 2\aaa_1      & 2\aaa_3 \\
         &             &              & 4\aaa_1
\end{array}\right)
\subset 4\aaa_1\amalg 2\aaa_3
$ & $18$   & $4_6^{+4}$ $\ast$ \\
\hline
\end{tabular}
\end{table}

\begin{table}
\begin{tabular}{|c||c|c|c|c|c|c|c|c|}
\hline
 {\bf n}& $G$   &$Deg$& $\rk S$ &$q_S$            \\
\hline
\hline
  $9$& $(C_2)^3$ &
$
\left(\begin{array}{cccc}
 2\aaa_1 & (4\aaa_1)_I & (4\aaa_1)_I  & 10\aaa_1 \\
         & 2\aaa_1     & (4\aaa_1)_I  & 10\aaa_1 \\
         &             & 2\aaa_1      & 10\aaa_1 \\
         &             &              & 8\aaa_1
\end{array}\right)
\subset 14\aaa_1
$ & $18$   & $4_6^{+4}$ \\
\hline
    &       &
$((2\aaa_1,2\aaa_1)_I,4\aaa_1,4\aaa_1)\subset 12\aaa_1$ & $18$   &
$4_6^{+4}$ $\ast$ \\
\hline
    &       &
$((2\aaa_1,2\aaa_1)_{II},4\aaa_1,4\aaa_1)\subset 12\aaa_1$ & $18$   &
$2_{II}^{-2},4_2^{-2}$ $\ast$ \\
\hline
     &       &
$
\left(\begin{array}{cccc}
 2\aaa_1 & 2\aaa_3     & (4\aaa_1)_I & 6\aaa_1 \\
         & 4\aaa_1     & 6\aaa_1  & 8\aaa_1 \\
         &             & 2\aaa_1      & 6\aaa_1 \\
         &             &              & 4\aaa_1
\end{array}\right)
\subset 2\aaa_3\amalg 6\aaa_1
$ & $18$   & $2_{II}^{-2},4_1^{+1},8_5^{-1}$ $\ast$ \\
\hline
     &       &
$
\left(\begin{array}{cccc}
 2\aaa_1 & 2\aaa_3     & (4\aaa_1)_I & 6\aaa_1 \\
         & 4\aaa_1     & 6\aaa_1     & 8\aaa_1 \\
         &             & 2\aaa_1     & 2\aaa_3 \\
         &             &              & 4\aaa_1
\end{array}\right)
\subset 4\aaa_3
$ & $18$   & $2_{II}^{-2},4_2^{-2}$ $\ast$ \\
\hline
    &       &
$((2\aaa_1,2\aaa_1)_I,4\aaa_1,8\aaa_1)
\subset 4\aaa_1\amalg 4\aaa_3$ & $18$   &
$2_{II}^{-2},4_2^{-2}$ $\ast$ \\
\hline
    &       &
$(2\aaa_1,4\aaa_1,4\aaa_1,4\aaa_1)\subset 14\aaa_1$ & $18$   &
$2_{II}^{+2},4_5^{-1},8_5^{-1}$ $\ast$ \\
\hline
    &       &
$(2\aaa_1,4\aaa_1,4\aaa_1,4\aaa_1)\subset 2\aaa_3\amalg 8\aaa_1$ & $18$   &
$2_{II}^{+2},4_6^{+2}$ $\ast$ \\
\hline
    &       &
$(4\aaa_1,4\aaa_1,4\aaa_1,4\aaa_1)\subset 16\aaa_1$ & $18$   &
$2_{II}^{+2},4_6^{+2}$ $\ast$ \\
\hline
     &       &
$
\left(\begin{array}{ccccc}
 2\aaa_1 & (4\aaa_1)_I & (4\aaa_1)_I  & (4\aaa_1)_I &  10\aaa_1 \\
         & 2\aaa_1     & (4\aaa_1)_I  & (4\aaa_1)_I &  10\aaa_1 \\
         &             & 2\aaa_1      & (4\aaa_1)_I &  10\aaa_1 \\
         &             &              & 2\aaa_1     &  10\aaa_1\\
         &             &              &             &   8\aaa_1
\end{array}\right)
\subset 16\aaa_1
$ & $19$   & $4_5^{+3}$ $\ast$ \\
\hline
     &       &
$
\left(\begin{array}{ccccc}
 2\aaa_1 & (4\aaa_1)_I & 6\aaa_1    & 6\aaa_1     &  6\aaa_1 \\
         & 2\aaa_1     & 6\aaa_1  & (4\aaa_1)_I &  6\aaa_1 \\
         &             & 4\aaa_1      & 6\aaa_1   &  8\aaa_1 \\
         &             &              & 2\aaa_1     & 2\aaa_3\\
         &             &              &             & 4\aaa_1
\end{array}\right)
\subset 8\aaa_1\amalg 2\aaa_3
$ & $19$   & $4_4^{-2},8_5^{-1}$ $\ast$ \\
\hline
     &       &
$
\left(\begin{array}{ccccc}
 2\aaa_1 & (4\aaa_1)_I & 6\aaa_1    & (4\aaa_1)_I     &  6\aaa_1 \\
         & 2\aaa_1     & 2\aaa_3  & (4\aaa_1)_I       &  6\aaa_1 \\
         &             & 4\aaa_1      & 6\aaa_1   &  8\aaa_1 \\
         &             &              & 2\aaa_1     & 2\aaa_3\\
         &             &              &             & 4\aaa_1
\end{array}\right)
\subset 2\aaa_1\amalg 4\aaa_3
$ & $19$   & $4_5^{+3}$ $\ast$ \\
\hline

\end{tabular}
\end{table}

\newpage

\begin{table}
\begin{tabular}{|c||c|c|c|c|c|c|c|c|}
\hline
 {\bf n}& $G$   &$Deg$& $\rk S$ &$q_S$            \\
\hline
\hline
  $9$& $(C_2)^3$ &
$
\left(\begin{array}{ccccc}
 2\aaa_1 & (4\aaa_1)_I & (4\aaa_1)_I    & 6\aaa_1     &  10\aaa_1 \\
         & 2\aaa_1     & (4\aaa_1)_I    & 6\aaa_1     &  10\aaa_1 \\
         &             & 2\aaa_1        & 6\aaa_1     &  10\aaa_1 \\
         &             &                & 4\aaa_1     &   4\aaa_3\\
         &             &              &               &   8\aaa_1
\end{array}\right)
\subset 6\aaa_1\amalg 4\aaa_3
$ & $19$   & $4_5^{+3}$ $\ast$ \\
\hline
      &       &
$((2\aaa_1,2\aaa_1)_{II},4\aaa_1,4\aaa_1,4\aaa_1)\subset 16\aaa_1$ & $19$   &
$2_{II}^{-2},8_5^{-1}$ $\ast$ \\
\hline
     &       &
$
\left(\begin{array}{ccccc}
 2\aaa_1 & 2\aaa_3     & (4\aaa_1)_I  &  6\aaa_1     &  6\aaa_1 \\
         & 4\aaa_1     & 6\aaa_1     &  8\aaa_1     &  8\aaa_1 \\
         &             & 2\aaa_1      & 2\aaa_3     &  6\aaa_1 \\
         &             &              & 4\aaa_1     &  8\aaa_1\\
         &             &              &             &  4\aaa_1
\end{array}\right)
\subset 4\aaa_3\amalg  4\aaa_1
$ & $19$   & $2_{II}^{-2},8_5^{-1}$ $\ast$ \\
\hline

\end{tabular}
\end{table}

\newpage

%%%%%%%%%%%%%%%%%%%%%%%%%%%%%%%%%%%%%%%%%%%%%%%%%%%
%%%%%%%%%%%%%%%%%%%%%%%%%%%%%%%%%%%%%%%%%%%%%%%%%%%%
%%%%%%%%%%%%%%%%%%%%%%%%%%%%%%%%%%%%%%%%%%%%%%%%%%%%

\newpage

\begin{table}
\label{table7}
\caption{Markings $S\subset N_j$ by Niemeier lattices, and lattices $S^\perp_{N_j}$
of degenerations of
codimension $\ge 2$ of K\"ahlerian K3 surfaces with
symplectic automorphism group $(C_2)^3$
in notations \cite{Nik8}, \cite{Nik9}, \cite{Nik10}.}
\end{table}

\medskip

{\bf n}=9, degeneration
$((2\aaa_1,2\aaa_1)\subset 4\aaa_1)_I$:

\begin{center}
\begin{tabular}{|c|c|c|c|c|c|c|c|c|}
\hline
$j$ & $21$ & $21$ & $21$ & $23$ & $23$ $\ast$ \\
\hline
$H$ & $H_{9,1}$ & $H_{9,2}$ & $H_{9,2}$ & $H_{9,3}$ & $H_{9,4}$ \\
\hline
orbits of & $(\alpha_{1,1},\alpha_{1,2}),...$ &
$(\alpha_{2,2},\alpha_{2,3}),...$ &
$(\alpha_{2,2},\alpha_{1,1}),...$&
$(\alpha_{6},\alpha_{8}),...$ &
$(\alpha_{1},\alpha_{8}),...$ \\
\hline
$(S^\perp_{N_j})^{(2)}$ & $4A_1$ & $4A_1$ & $2A_1$ &$2A_1$ & $\{0\}$  \\
\hline
\end{tabular}
\end{center}

\medskip

{\bf n}=9, degeneration
$((2\aaa_1,2\aaa_1)\subset 4\aaa_1)_{II}$:

\begin{center}
\begin{tabular}{|c|c|c|c|c|c|c|c|c|}
\hline
$j$ & $23$ $\ast$ \\
\hline
$H$ & $H_{9,1}$ \\
\hline
orbits of & $(\alpha_{5},\alpha_{6})$ \\
\hline
$(S^\perp_{N_j})^{(2)}$ & $4A_1$  \\
\hline
\end{tabular}
\end{center}

\medskip

{\bf n}=9, degeneration
$(2\aaa_1,4\aaa_1)\subset 6\aaa_1$:

\begin{center}
\begin{tabular}{|c|c|c|c|c|c|c|c|c|}
\hline
$j$ & $21$ & $21$ & $21$ & $23$ & $23$ $\ast$ \\
\hline
$H$ & $H_{9,1}$ & $H_{9,2}$ & $H_{9,2}$ & $H_{9,1}$ & $H_{9,4}$ \\
\hline
orbits of & $(\alpha_{1,1},\alpha_{2,3}),...$ &
$(\alpha_{2,3},\alpha_{1,2}),...$ &
$(\alpha_{1,1},\alpha_{1,2}),...$&
$(\alpha_{5},\alpha_{2}),...$ &
$(\alpha_{1},\alpha_{2}),...$ \\
\hline
$(S^\perp_{N_j})^{(2)}$ & $6A_1$ & $4A_1$ & $2A_1$ &$4A_1$ & $\{0\}$  \\
\hline
\end{tabular}
\end{center}

\medskip

{\bf n}=9, degeneration
$(2\aaa_1,4\aaa_1)\subset 2\aaa_3$:

\begin{center}
\begin{tabular}{|c|c|c|c|c|c|c|c|c|}
\hline
$j$ & $21$ $\ast$ \\
\hline
$H$ & $H_{9,2}$ \\
\hline
orbits of & $(\alpha_{2,2},\alpha_{1,2}),...$ \\
\hline
$(S^\perp_{N_j})^{(2)}$ & $4A_1$  \\
\hline
\end{tabular}
\end{center}

\medskip

{\bf n}=9, degeneration
$(2\aaa_1,8\aaa_1)\subset 10\aaa_1$:

\begin{center}
\begin{tabular}{|c|c|c|c|c|c|c|c|c|}
\hline
$j$ & $21$ & $23$ $\ast$ \\
\hline
$H$ & $H_{9,1}$ & $H_{9,3}$  \\
\hline
orbits of & $(\alpha_{1,1},\alpha_{1,3}),...$ &
$(\alpha_{6},\alpha_{2}),...$ \\
\hline
$(S^\perp_{N_j})^{(2)}$ & $6A_1$ & $2A_1$   \\
\hline
\end{tabular}
\end{center}

\medskip

{\bf n}=9, degeneration
$(4\aaa_1,4\aaa_1)\subset 8\aaa_1$:

\begin{center}
\begin{tabular}{|c|c|c|c|c|c|c|c|c|}
\hline
$j$ & $21$ & $23$ \\
\hline
$H$ & $H_{9,2}$ & $H_{9,1}$  \\
\hline
orbits of & $(\alpha_{1,2},\alpha_{1,3}),...$ &
$(\alpha_{2},\alpha_{3}),...$ \\
\hline
$(S^\perp_{N_j})^{(2)}$ & $4A_1$ & $4A_1$   \\
\hline
\end{tabular}
\end{center}

\medskip

{\bf n}=9, degeneration
$(4\aaa_1,8\aaa_1)\subset 4\aaa_3$:

\begin{center}
\begin{tabular}{|c|c|c|c|c|c|c|c|c|}
\hline
$j$ & $21$ $\ast$ \\
\hline
$H$ & $H_{9,1}$ \\
\hline
orbits of & $(\alpha_{2,3},\alpha_{1,3})$ \\
\hline
$(S^\perp_{N_j})^{(2)}$ & $8A_1$  \\
\hline
\end{tabular}
\end{center}

\medskip

{\bf n}=9, degeneration
$(8\aaa_1,8\aaa_1)\subset 16\aaa_1$:

\begin{center}
\begin{tabular}{|c|c|c|c|c|c|c|c|c|}
\hline
$j$ & $23$ $\ast$ \\
\hline
$H$ & $H_{9,2}$ \\
\hline
orbits of & $(\alpha_{2},\alpha_{4})$ \\
\hline
$(S^\perp_{N_j})^{(2)}$ & $8A_1$  \\
\hline
\end{tabular}
\end{center}

\medskip

{\bf n}=9, degeneration
$
\left(\begin{array}{ccc}
 2\aaa_1 & (4\aaa_1)_I & (4\aaa_1)_I  \\
         & 2\aaa_1 & (4\aaa_1)_I   \\
         &         & 2\aaa_1
\end{array}\right)
\subset 6\aaa_1
$:

\begin{center}
\begin{tabular}{|c|c|c|c|c|c|c|c|c|}
\hline
$j$ & $21$ & $21$ & $23$ & $23$ $\ast$ \\
\hline
$H$ & $H_{9,1}$ & $H_{9,2}$ & $H_{9,3}$ & $H_{9,4}$  \\
\hline
orbits of & $(\alpha_{1,1},\alpha_{1,2},\alpha_{1,4}),...$ &
$(\alpha_{2,2},\alpha_{2,3},\alpha_{1,1}),...$ &
$(\alpha_{6},\alpha_{8},\alpha_{10}),...$ &
$(\alpha_{1},\alpha_{8},\alpha_{10}),... $ \\
\hline
$(S^\perp_{N_j})^{(2)}$ & $2A_1$ & $2A_1$ &$2A_1$ & $\{0\}$  \\
\hline
\end{tabular}
\end{center}

\medskip

{\bf n}=9, degeneration
$((2\aaa_1,2\aaa_1)_I,4\aaa_1)\subset 8\aaa_1$:

\begin{center}
\begin{tabular}{|c|c|c|c|c|c|c|c|c|}
\hline
$j$ & $21$ & $21$ & $21$ & $23$ $\ast$ \\
\hline
$H$ & $H_{9,1}$ & $H_{9,2}$ & $H_{9,2}$ & $H_{9,4}$  \\
\hline
orbits of & $(\alpha_{1,1},\alpha_{1,2},\alpha_{2,3}),...$ &
$(\alpha_{2,3},\alpha_{2,7},\alpha_{1,2}),...$ &
$(\alpha_{2,3},\alpha_{1,1},\alpha_{1,2}),...$ &
$(\alpha_{1},\alpha_{8},\alpha_{2}),...$ \\
\hline
$(S^\perp_{N_j})^{(2)}$ & $4A_1$ & $4A_1$ &$2A_1$ & $\{0\}$  \\
\hline
\end{tabular}
\end{center}

\medskip

{\bf n}=9, degeneration
$((2\aaa_1,2\aaa_1)_{II},4\aaa_1)\subset 8\aaa_1$:

\begin{center}
\begin{tabular}{|c|c|c|c|c|c|c|c|c|}
\hline
$j$ & $23$ $\ast$ \\
\hline
$H$ & $H_{9,1}$  \\
\hline
orbits of & $(\alpha_{5},\alpha_{6},\alpha_{2}),...$ \\
\hline
$(S^\perp_{N_j})^{(2)}$ & $4A_1$   \\
\hline
\end{tabular}
\end{center}

\medskip

{\bf n}=9, degeneration
$
\left(\begin{array}{ccc}
 2\aaa_1 & (4\aaa_1)_I & 6\aaa_1  \\
         & 2\aaa_1 & 2\aaa_3   \\
         &         & 4\aaa_1
\end{array}\right)
\subset 2\aaa_1\amalg 2\aaa_3
$:

\begin{center}
\begin{tabular}{|c|c|c|c|c|c|c|c|c|}
\hline
$j$ & $21$ $\ast$ & $21$ \\
\hline
$H$ & $H_{9,2}$ & $H_{9,2}$   \\
\hline
orbits of & $(\alpha_{2,3},\alpha_{2,2},\alpha_{1,2}),...$ &
$(\alpha_{1,1},\alpha_{2,2},\alpha_{1,2}),...$ \\
\hline
$(S^\perp_{N_j})^{(2)}$ & $4A_1$ & $2A_1$  \\
\hline
\end{tabular}
\end{center}

\medskip

{\bf n}=9, degeneration
$((2\aaa_1,2\aaa_1)_I,8\aaa_1)\subset 12\aaa_1$:

\begin{center}
\begin{tabular}{|c|c|c|c|c|c|c|c|c|}
\hline
$j$ & $21$ & $23$ $\ast$ \\
\hline
$H$ & $H_{9,1}$ & $H_{9,3}$   \\
\hline
orbits of & $(\alpha_{1,1},\alpha_{1,2},\alpha_{1,3}),...$ &
$(\alpha_{6},\alpha_{8},\alpha_{2}),...$ \\
\hline
$(S^\perp_{N_j})^{(2)}$ & $4A_1$ & $2A_1$  \\
\hline
\end{tabular}
\end{center}

\medskip

{\bf n}=9, degeneration
$(2\aaa_1,4\aaa_1,4\aaa_1)\subset 10\aaa_1$:

\begin{center}
\begin{tabular}{|c|c|c|c|c|c|c|c|c|}
\hline
$j$ & $21$ & $21$ $\ast$ & $23$ \\
\hline
$H$ & $H_{9,2}$ & $H_{9,2}$ & $H_{9,1}$  \\
\hline
orbits of & $(\alpha_{2,7},\alpha_{1,2},\alpha_{1,3}),...$ &
$(\alpha_{1,1},\alpha_{1,2},\alpha_{1,3}),...$ &
$(\alpha_{5},\alpha_{2},\alpha_{3}),...$ \\
\hline
$(S^\perp_{N_j})^{(2)}$ & $4A_1$ & $2A_1$ &$4A_1$   \\
\hline
\end{tabular}
\end{center}

\medskip

{\bf n}=9, degeneration
$(2\aaa_1,4\aaa_1,4\aaa_1)\subset 2\aaa_3\amalg 4\aaa_1$:

\begin{center}
\begin{tabular}{|c|c|c|c|c|c|c|c|c|}
\hline
$j$ & $21$ $\ast$ \\
\hline
$H$ & $H_{9,2}$  \\
\hline
orbits of & $(\alpha_{2,2},\alpha_{1,2},\alpha_{1,3}),...$ \\
\hline
$(S^\perp_{N_j})^{(2)}$ & $4A_1$   \\
\hline
\end{tabular}
\end{center}

\medskip

{\bf n}=9, degeneration
$(2\aaa_1,4\aaa_1,8\aaa_1)\subset 2\aaa_1\amalg 4\aaa_3$:

\begin{center}
\begin{tabular}{|c|c|c|c|c|c|c|c|c|}
\hline
$j$ & $21$ $\ast$ \\
\hline
$H$ & $H_{9,1}$  \\
\hline
orbits of & $(\alpha_{1,1},\alpha_{2,3},\alpha_{1,3}),...$ \\
\hline
$(S^\perp_{N_j})^{(2)}$ & $6A_1$   \\
\hline
\end{tabular}
\end{center}

\medskip

{\bf n}=9, degeneration
$(4\aaa_1,4\aaa_1,4\aaa_1)\subset 12\aaa_1$: 

\begin{center}
\begin{tabular}{|c|c|c|c|c|c|c|c|c|}
\hline
$j$ & $21$ & $23$ \\
\hline
$H$ & $H_{9,2}$ & $H_{9,1}$   \\
\hline
orbits of & $(\alpha_{1,2},\alpha_{1,3},\alpha_{1,7}),...$ &
$(\alpha_{2},\alpha_{3},\alpha_{4}),...$ \\
\hline
$(S^\perp_{N_j})^{(2)}$ & $4A_1$ & $4A_1$  \\
\hline
\end{tabular}
\end{center}

\medskip

{\bf n}=9, degeneration
$
\left(\begin{array}{cccc}
 2\aaa_1 & (4\aaa_1)_I & (4\aaa_1)_I  & (4\aaa_1)_I \\
         & 2\aaa_1 & (4\aaa_1)_I & (4\aaa_1)_I  \\
         &         & 2\aaa_1 &    (4\aaa_1)_I \\
         &         &         & 2\aaa_1
\end{array}\right)
\subset 8\aaa_1
$:

\begin{center}
\begin{tabular}{|c|c|c|c|c|c|c|c|c|}
\hline
$j$ & $23$ $\ast$ \\
\hline
$H$ & $H_{9,3}$  \\
\hline
orbits of & $(\alpha_{6},\alpha_{8},\alpha_{10},\alpha_{16}),...$ \\
\hline
$(S^\perp_{N_j})^{(2)}$ & $2A_1$   \\
\hline
\end{tabular}
\end{center}

\medskip

{\bf n}=9, degeneration
$
\left(\begin{array}{cccc}
 2\aaa_1 & (4\aaa_1)_I & (4\aaa_1)_I  & 6\aaa_1 \\
         & 2\aaa_1     & (4\aaa_1)_I  & 6\aaa_1 \\
         &             & 2\aaa_1      & 6\aaa_1 \\
         &             &              & 4\aaa_1
\end{array}\right)
\subset 10\aaa_1
$:

\begin{center}
\begin{tabular}{|c|c|c|c|c|c|c|c|c|}
\hline
$j$ & $21$ & $21$  & $23$ $\ast$ \\
\hline
$H$ & $H_{9,1}$ & $H_{9,2}$ & $H_{9,4}$  \\
\hline
orbits of & $(\alpha_{1,1},\alpha_{1,2},\alpha_{1,4},\alpha_{2,3}),...$ &
$(\alpha_{2,3},\alpha_{2,7},\alpha_{1,1},\alpha_{1,2}),...$ &
$(\alpha_{1},\alpha_{8},\alpha_{10},\alpha_2 ),...$ \\
\hline
$(S^\perp_{N_j})^{(2)}$ & $2A_1$ & $2A_1$ &$\{0\}$   \\
\hline
\end{tabular}
\end{center}

\medskip

{\bf n}=9, degeneration
$
\left(\begin{array}{cccc}
 2\aaa_1 & (4\aaa_1)_I & (4\aaa_1)_I  & 6\aaa_1 \\
         & 2\aaa_1     & (4\aaa_1)_I  & 6\aaa_1 \\
         &             & 2\aaa_1      & 2\aaa_3 \\
         &             &              & 4\aaa_1
\end{array}\right)
\subset 4\aaa_1\amalg 2\aaa_3
$:

\begin{center}
\begin{tabular}{|c|c|c|c|c|c|c|c|c|}
\hline
$j$ & $21$ $\ast$ \\
\hline
$H$ & $H_{9,2}$  \\
\hline
orbits of & $(\alpha_{1,1},\alpha_{2,3},\alpha_{2,2},\alpha_{1,2}),...$ \\
\hline
$(S^\perp_{N_j})^{(2)}$ & $2A_1$   \\
\hline
\end{tabular}
\end{center}

\medskip

{\bf n}=9, degeneration
$
\left(\begin{array}{cccc}
 2\aaa_1 & (4\aaa_1)_I & (4\aaa_1)_I  & 10\aaa_1 \\
         & 2\aaa_1     & (4\aaa_1)_I  & 10\aaa_1 \\
         &             & 2\aaa_1      & 10\aaa_1 \\
         &             &              & 8\aaa_1
\end{array}\right)
\subset 14\aaa_1
$:

\begin{center}
\begin{tabular}{|c|c|c|c|c|c|c|c|c|}
\hline
$j$ & $21$  & $23$ \\
\hline
$H$ & $H_{9,1}$ & $H_{9,3}$  \\
\hline
orbits of & $(\alpha_{1,1},\alpha_{1,2},\alpha_{1,4},\alpha_{1,3}),...$ &
$(\alpha_{6},\alpha_{8},\alpha_{10},\alpha_2 ),...$ \\
\hline
$(S^\perp_{N_j})^{(2)}$ & $2A_1$ & $2A_1$   \\
\hline
\end{tabular}
\end{center}

\medskip

{\bf n}=9, degeneration
$((2\aaa_1,2\aaa_1)_I,4\aaa_1,4\aaa_1)\subset 12\aaa_1$:

\begin{center}
\begin{tabular}{|c|c|c|c|c|c|c|c|c|}
\hline
$j$ & $21$ $\ast$ \\
\hline
$H$ & $H_{9,2}$  \\
\hline
orbits of & $(\alpha_{2,7},\alpha_{1,1},\alpha_{1,2},\alpha_{1,3}),...$ \\
\hline
$(S^\perp_{N_j})^{(2)}$ & $2A_1$   \\
\hline
\end{tabular}
\end{center}

\medskip

{\bf n}=9, degeneration
$((2\aaa_1,2\aaa_1)_{II},4\aaa_1,4\aaa_1)\subset 12\aaa_1$:

\begin{center}
\begin{tabular}{|c|c|c|c|c|c|c|c|c|}
\hline
$j$ & $23$ $\ast$ \\
\hline
$H$ & $H_{9,1}$  \\
\hline
orbits of & $(\alpha_{5},\alpha_{6},\alpha_{2},\alpha_{3}),...$ \\
\hline
$(S^\perp_{N_j})^{(2)}$ & $4A_1$   \\
\hline
\end{tabular}
\end{center}

\medskip

{\bf n}=9, degeneration
$
\left(\begin{array}{cccc}
 2\aaa_1 & 2\aaa_3     & (4\aaa_1)_I & 6\aaa_1 \\
         & 4\aaa_1     & 6\aaa_1  & 8\aaa_1 \\
         &             & 2\aaa_1      & 6\aaa_1 \\
         &             &              & 4\aaa_1
\end{array}\right)
\subset 2\aaa_3\amalg 6\aaa_1
$:

\begin{center}
\begin{tabular}{|c|c|c|c|c|c|c|c|c|}
\hline
$j$ & $21$ $\ast$  & $21$ \\
\hline
$H$ & $H_{9,2}$ & $H_{9,2}$  \\
\hline
orbits of & $(\alpha_{2,2},\alpha_{1,2},\alpha_{2,7},\alpha_{1,3}),...$ &
$(\alpha_{2,2},\alpha_{1,2},\alpha_{1,1},\alpha_{1,3}),...$ \\
\hline
$(S^\perp_{N_j})^{(2)}$ & $4A_1$ & $2A_1$   \\
\hline
\end{tabular}
\end{center}

\medskip

{\bf n}=9, degeneration
$
\left(\begin{array}{cccc}
 2\aaa_1 & 2\aaa_3     & (4\aaa_1)_I & 6\aaa_1 \\
         & 4\aaa_1     & 6\aaa_1     & 8\aaa_1 \\
         &             & 2\aaa_1     & 2\aaa_3 \\
         &             &              & 4\aaa_1
\end{array}\right)
\subset 4\aaa_3
$:

\begin{center}
\begin{tabular}{|c|c|c|c|c|c|c|c|c|}
\hline
$j$ & $21$ $\ast$ \\
\hline
$H$ & $H_{9,2}$  \\
\hline
orbits of & $(\alpha_{2,2},\alpha_{1,2},\alpha_{2,3},\alpha_{1,3}),...$ \\
\hline
$(S^\perp_{N_j})^{(2)}$ & $4A_1$   \\
\hline
\end{tabular}
\end{center}

\medskip

{\bf n}=9, degeneration
$((2\aaa_1,2\aaa_1)_I,4\aaa_1,8\aaa_1)
\subset 4\aaa_1\amalg 4\aaa_3$:

\begin{center}
\begin{tabular}{|c|c|c|c|c|c|c|c|c|}
\hline
$j$ & $21$ $\ast$ \\
\hline
$H$ & $H_{9,1}$  \\
\hline
orbits of & $(\alpha_{1,1},\alpha_{1,2},\alpha_{2,3},\alpha_{1,3}),...$ \\
\hline
$(S^\perp_{N_j})^{(2)}$ & $4A_1$   \\
\hline
\end{tabular}
\end{center}

\medskip

{\bf n}=9, degeneration
$(2\aaa_1,4\aaa_1,4\aaa_1,4\aaa_1)\subset 14\aaa_1$:

\begin{center}
\begin{tabular}{|c|c|c|c|c|c|c|c|c|}
\hline
$j$ & $21$  & $23$ $\ast$\\
\hline
$H$ & $H_{9,2}$ & $H_{9,1}$  \\
\hline
orbits of & $(\alpha_{1,1},\alpha_{1,2},\alpha_{1,3},\alpha_{1,7}),...$ &
$(\alpha_{5},\alpha_{2},\alpha_{3},\alpha_4),...$ \\
\hline
$(S^\perp_{N_j})^{(2)}$ & $2A_1$ & $4A_1$   \\
\hline
\end{tabular}
\end{center}

\medskip

{\bf n}=9, degeneration
$(2\aaa_1,4\aaa_1,4\aaa_1,4\aaa_1)\subset 2\aaa_3\amalg 8\aaa_1$:

\begin{center}
\begin{tabular}{|c|c|c|c|c|c|c|c|c|}
\hline
$j$ & $21$ $\ast$ \\
\hline
$H$ & $H_{9,2}$  \\
\hline
orbits of & $(\alpha_{2,2},\alpha_{1,2},\alpha_{1,3},\alpha_{1,7}),...$ \\
\hline
$(S^\perp_{N_j})^{(2)}$ & $4A_1$   \\
\hline
\end{tabular}
\end{center}

\medskip

{\bf n}=9, degeneration
$(4\aaa_1,4\aaa_1,4\aaa_1,4\aaa_1)\subset 16\aaa_1$:

\begin{center}
\begin{tabular}{|c|c|c|c|c|c|c|c|c|}
\hline
$j$ & $23$ $\ast$ \\
\hline
$H$ & $H_{9,1}$  \\
\hline
orbits of & $(\alpha_{2},\alpha_{3},\alpha_{4},\alpha_{7}),...$ \\
\hline
$(S^\perp_{N_j})^{(2)}$ & $4A_1$   \\
\hline
\end{tabular}
\end{center}

\medskip

{\bf n}=9, degeneration
$
\left(\begin{array}{ccccc}
 2\aaa_1 & (4\aaa_1)_I & (4\aaa_1)_I  & (4\aaa_1)_I &  10\aaa_1 \\
         & 2\aaa_1     & (4\aaa_1)_I  & (4\aaa_1)_I &  10\aaa_1 \\
         &             & 2\aaa_1      & (4\aaa_1)_I &  10\aaa_1 \\
         &             &              & 2\aaa_1     &  10\aaa_1\\
         &             &              &             &   8\aaa_1
\end{array}\right)
\subset 16\aaa_1
$:

\begin{center}
\begin{tabular}{|c|c|c|c|c|c|c|c|c|}
\hline
$j$ & $23$ $\ast$ \\
\hline
$H$ & $H_{9,3}$  \\
\hline
orbits of & $(\alpha_{6},\alpha_{8},\alpha_{10},\alpha_{16},\alpha_{2}),...$ \\
\hline
$(S^\perp_{N_j})^{(2)}$ & $2A_1$   \\
\hline
\end{tabular}
\end{center}

\medskip

{\bf n}=9, degeneration
$
\left(\begin{array}{ccccc}
 2\aaa_1 & (4\aaa_1)_I & 6\aaa_1    & 6\aaa_1     &  6\aaa_1 \\
         & 2\aaa_1     & 6\aaa_1  & (4\aaa_1)_I &  6\aaa_1 \\
         &             & 4\aaa_1      & 6\aaa_1   &  8\aaa_1 \\
         &             &              & 2\aaa_1     & 2\aaa_3\\
         &             &              &             & 4\aaa_1
\end{array}\right)
\subset 8\aaa_1\amalg 2\aaa_3
$:

\begin{center}
\begin{tabular}{|c|c|c|c|c|c|c|c|c|}
\hline
$j$ & $21$ $\ast$ \\
\hline
$H$ & $H_{9,2}$  \\
\hline
orbits of & $(\alpha_{1,1},\alpha_{2,7},\alpha_{1,3},\alpha_{2,2},\alpha_{1,2}),...$ \\
\hline
$(S^\perp_{N_j})^{(2)}$ & $2A_1$   \\
\hline
\end{tabular}
\end{center}

\medskip

{\bf n}=9, degeneration
$
\left(\begin{array}{ccccc}
 2\aaa_1 & (4\aaa_1)_I & 6\aaa_1    & (4\aaa_1)_I     &  6\aaa_1 \\
         & 2\aaa_1     & 2\aaa_3  & (4\aaa_1)_I       &  6\aaa_1 \\
         &             & 4\aaa_1      & 6\aaa_1   &  8\aaa_1 \\
         &             &              & 2\aaa_1     & 2\aaa_3\\
         &             &              &             & 4\aaa_1
\end{array}\right)
\subset 2\aaa_1\amalg 4\aaa_3
$:

\begin{center}
\begin{tabular}{|c|c|c|c|c|c|c|c|c|}
\hline
$j$ & $21$ $\ast$ \\
\hline
$H$ & $H_{9,2}$  \\
\hline
orbits of & $(\alpha_{1,1},\alpha_{2,2},\alpha_{1,2},\alpha_{2,3},\alpha_{1,3}),...$ \\
\hline
$(S^\perp_{N_j})^{(2)}$ & $2A_1$   \\
\hline
\end{tabular}
\end{center}

\medskip

{\bf n}=9, degeneration
$
\left(\begin{array}{ccccc}
 2\aaa_1 & (4\aaa_1)_I & (4\aaa_1)_I    & 6\aaa_1     &  10\aaa_1 \\
         & 2\aaa_1     & (4\aaa_1)_I    & 6\aaa_1     &  10\aaa_1 \\
         &             & 2\aaa_1        & 6\aaa_1     &  10\aaa_1 \\
         &             &                & 4\aaa_1     &   4\aaa_3\\
         &             &              &               &   8\aaa_1
\end{array}\right)
\subset 6\aaa_1\amalg 4\aaa_3
$:

\begin{center}
\begin{tabular}{|c|c|c|c|c|c|c|c|c|}
\hline
$j$ & $21$ $\ast$ \\
\hline
$H$ & $H_{9,1}$  \\
\hline
orbits of & $(\alpha_{1,1},\alpha_{1,2},\alpha_{1,4},\alpha_{2,3},\alpha_{1,3}),...$ \\
\hline
$(S^\perp_{N_j})^{(2)}$ & $2A_1$   \\
\hline
\end{tabular}
\end{center}

\medskip

{\bf n}=9, degeneration
$((2\aaa_1,2\aaa_1)_{II},4\aaa_1,4\aaa_1,4\aaa_1)\subset 16\aaa_1$:

\begin{center}
\begin{tabular}{|c|c|c|c|c|c|c|c|c|}
\hline
$j$ & $23$ $\ast$ \\
\hline
$H$ & $H_{9,1}$  \\
\hline
orbits of & $(\alpha_{5},\alpha_{6},\alpha_{2},\alpha_{3},\alpha_{4}),...$ \\
\hline
$(S^\perp_{N_j})^{(2)}$ & $4A_1$   \\
\hline
\end{tabular}
\end{center}

\medskip

{\bf n}=9, degeneration
$
\left(\begin{array}{ccccc}
 2\aaa_1 & 2\aaa_3     & (4\aaa_1)_I  &  6\aaa_1     &  6\aaa_1 \\
         & 4\aaa_1     & 6\aaa_1     &  8\aaa_1     &  8\aaa_1 \\
         &             & 2\aaa_1      & 2\aaa_3     &  6\aaa_1 \\
         &             &              & 4\aaa_1     &  8\aaa_1\\
         &             &              &             &  4\aaa_1
\end{array}\right)
\subset 4\aaa_3\amalg  4\aaa_1
$:

\begin{center}
\begin{tabular}{|c|c|c|c|c|c|c|c|c|}
\hline
$j$ & $21$ $\ast$ \\
\hline
$H$ & $H_{9,2}$  \\
\hline
orbits of & $(\alpha_{2,2},\alpha_{1,2},\alpha_{2,3},\alpha_{1,3},\alpha_{1,7}),...$ \\
\hline
$(S^\perp_{N_j})^{(2)}$ & $4A_1$   \\
\hline
\end{tabular}
\end{center}

\medskip

\section{Final remarks }
\label{sec5.remarks}

We hope to give more details and applications in further variants of 
the paper and further publications. 

We hope to consider remaining groups $D_6$, $C_4$, $(C_2)^2$, $C_3$, $C_2$ 
later as well.

V.V. Nikulin
\par Steklov Mathematical Institute,
\par ul. Gubkina 8, Moscow 117966, GSP-1, Russia;

\vskip5pt

Deptm. of Pure Mathem. The University of
Liverpool, Liverpool\par L69 3BX, UK
\par

\vskip5pt

nikulin@mi.ras.ru\, \ \
vnikulin@liv.ac.uk \, \ \  vvnikulin@list.ru
%%%%%%%%%%%%%%%%

Personal page: http://vnikulin.com

\end{document}